%% file: Almutari,_Gluck-_Infinitely_many_arXiv.tex
\documentclass[reqno]{amsart}

\input{structure.tex}

\newif\ifdetails

\title[Infinitely many solutions]{Infinitely many solutions to a conformally invariant elliptic equation with Choquard-type nonlinearity}
\date{\today}
\author{Mona Almutairi}
\address{Department of Mathematics, College of Science\\
Qassim University\\
	Burayday 51452, Saudi Arabia}
\email{mm.almutairi@qu.edu.sa}
\author{Mathew Gluck}
\thanks{The second author is supported by the National Science Foundation Grant No. DMS-2418889.}
\address{Southern Illinois University \\
	School of Mathematical and Statistical Sciences\\ 
	Carbondale, IL, U.S.A.
	}
\email{mathew.gluck@siu.edu}

\renewcommand{\div}{\text{div}}

\begin{document}
\begin{abstract}
The existence of an unbounded sequence of solutions to a conformally invariant elliptic equation having nonlocal critical-power nonlinearity is established. The primary obstacle to establishing existence of solutions is the failure of compactness in the Sobolev embedding. To overcome this obstacle, the problem under consideration is lifted to an equivalent problem on the standard sphere so that the symmetries of the sphere can be leveraged. Two classes of symmetries are considered and for each class of symmetries, an unbounded sequence of solutions to the lifted problem with the prescribed symmetries is produced. One class of symmetries always exists and the corresponding solutions are guaranteed to be sign-changing whenever a suitable relationship between the dimension and the nonlocality parameter holds. The other class of symmetries need not always exist but when it exists, the corresponding solutions are guaranteed to be sign-changing.
\end{abstract}
\maketitle
%

\section{Introduction}
This work is dedicated to the study of the nonlinear nonlocal problem
\begin{equation}
\label{eq:choquard}
\begin{cases}
	-\lap u  = I_\mu[ |u|^{2^*_\mu}]|u|^{2^*_\mu - 2}u \qquad \text{ in }\bb R^n\\
	u\in D^{1, 2}(\bb R^n)
\end{cases}
\end{equation}
where $n\geq 3$, $\mu\in (0, n)$, and 
\begin{equation}
\label{eq:the_exponent}
	2^*_\mu = \frac{2n - \mu}{n - 2}
\end{equation} 
is the critical exponent for problem \eqref{eq:choquard}. Here $I_\mu$ is the convolution operator defined by 
\begin{equation}
\label{eq:Reisz_potential_operator}
	I_\mu f(x) 
	= \int_{\bb R^n}\frac{f(y)}{|x - y|^\mu}\; \d y\qquad \text{for }x\in \bb R^n. 
\end{equation}
The nonlinearity appearing on the right-hand side of equation \eqref{eq:choquard} is commonly referred to as a nonlinearity of \emph{Choquard} type. For physical and historical motivations for studying problems with Choquard type nonlinearities we refer the reader to the survey \cite{Moroz2017} and the references therein. The motivation for the present work comes from the observation that problem \eqref{eq:choquard} has the same symmetries as the problem
\begin{equation}
\label{eq:mu=0_choquard}
\begin{cases}
	-\lap u = |u|^{2^*- 2}u \qquad \text{ in }\bb R^n\\
	u\in D^{1, 2}(\bb R^n),  
\end{cases}
\end{equation}
where in the notation of \eqref{eq:the_exponent}, 
\begin{equation}
\label{eq:2_star}
	2^* = 2^*_0 = \frac{2n}{n - 2}.
\end{equation} 
Problem \eqref{eq:mu=0_choquard} has a rich history due in part to the fact that its solutions serve as models for how compactness in the critical Sobolev embedding can fail. Moreover, up to constant multiple, problem \eqref{eq:mu=0_choquard} is formally obtained from problem \eqref{eq:choquard} both by letting $\mu\to 0^+$ and by letting $\mu\to n^-$. The classification of positive solutions to \eqref{eq:mu=0_choquard} has been established in \cite{GidasNiNirenberg1979, CaffarelliGidasSpruck1989} where it was shown that every positive solution to \eqref{eq:mu=0_choquard} is of the form
\begin{equation}
\label{eq:scaled_translated_bubble}
	U_{x_0, \lambda}(x) 
	= \lambda^{-(n - 2)/2}U\left(\frac{x - x_0}\lambda\right)
\end{equation}
for some $(x_0, \lambda)\in \bb R^n\times(0, \infty)$, where
\begin{equation}
\label{eq:UnscaledBubble}
	U(x) = c_n(1 + |x|^2)^{-(n - 2)/2} 
\end{equation}
and $c_n = (n(n -2))^{(n -2)/4}$. According to the folklore, after the release of the classification result in \cite{GidasNiNirenberg1979} many authors attempted to show that every solution to \eqref{eq:mu=0_choquard} that is positive somewhere is necessarily of the form \eqref{eq:scaled_translated_bubble}. The following theorem due to W. Ding implies that none of these attempts had potential for success.
\begin{oldtheorem}[\cite{Ding1986}]
\label{oldtheorem:Ding}
Let $n \geq 3$. There is a sequence $(u_i)_{i = 1}^\infty$ of solutions to \eqref{eq:mu=0_choquard} for which $\norm{\Grad u_i}_{L^2(\bb R^n)}\to\infty$ as $i\to \infty$. 
\end{oldtheorem}
\noindent Because every positive solution to \eqref{eq:mu=0_choquard} is of the form \eqref{eq:scaled_translated_bubble} for some $(x_0, \lambda)\in \bb R^n\times (0,\infty)$ and because the value of $\norm{\Grad U_{x_0,\lambda}}_{L^2(\bb R^n)}$ is independent of $x_0$ and $\lambda$, all but finitely many of the functions whose existence is guaranteed by Theorem \ref{oldtheorem:Ding} are sign-changing. 

Since the release of Theorem \ref{oldtheorem:Ding}, various works have been dedicated to extending Theorem \ref{oldtheorem:Ding} to different settings. One collection of such works concerns the family of problems 
\begin{equation}
\label{eq:mu=0_s_choquard}
	(-\lap)^s u = |u|^{\frac 4{n - 2s}}u
	\qquad \text{ in }\bb R^n, 
\end{equation}
indexed by $s\in (0, \frac n 2)$ (problem \eqref{eq:mu=0_choquard} corresponds to $s = 1$). The classification of positive solutions to problem \eqref{eq:mu=0_s_choquard} for the full range of $s\in (0, \frac n 2)$ has been obtained through the combined works of \cite{GidasNiNirenberg1979, CaffarelliGidasSpruck1989,Lin1998, Xu2000, WeiXu1999, ChenLiOu2006}. For  $s\in \{ 1, 2, \ldots, \lfloor\frac n 2\rfloor\}$, the existence of an unbounded sequence of solutions to \eqref{eq:mu=0_s_choquard} was established by Bartsch, Schneider and Weth \cite{BartschSchneiderWeth2004}. 
This result was extended by Abreu, Barbosa and Ramirez in \cite{AbreuEtAl2021} to include the full range of $s \in (0, n/2)$.%

Aside from extending Theorem \ref{oldtheorem:Ding} from the case $s= 1$ to the full range $s\in (0, n/2)$, a variety of similar problems in different contexts have been considered. For example, in \cite{Gluck2022} the existence of a sequence of solutions $(f_i)_{i = 1}^\infty\subset L^{\frac{2n}{2n - \mu}}(\bb R^n)$ to the problem
\begin{equation}
\label{eq:HLS_euler_lagrange}
	|f(y)|^{-\frac{2(n - \mu)}{2n - \mu}}f(y)
	= \int_{\bb R^n}\frac{|I_\mu f(x)|^{\frac{2(n- \mu)}{\mu}}I_\mu f(x)}{|x - y|^\mu}\; \d x
	\qquad y\in \bb R^n  
\end{equation}
for which $\|f_i\|_{L^{\frac{2n}{2n -\mu}}(\bb R^n)}\to \infty$ was established. Here $\mu\in (0, n)$ and $I_\mu$ is as in equation \eqref{eq:Reisz_potential_operator}. Just as the problems \eqref{eq:mu=0_s_choquard} arise as Euler-Lagrange equations for the critical points of suitable Sobolev quotients, problem \eqref{eq:HLS_euler_lagrange} arises as the Euler-Lagrange equation for the critical points of the diagonal Hardy-Littlewood-Sobolev quotient $f\mapsto \|f\|_{2n/(2n - \mu)}^{-1}\|I_\mu f\|_{2n/\mu}$. Other works in the spirit of Ding's Theorem include an analog of Ding's result for the spinorial Yamabe problem on the standard sphere \cite{Maalaoui2016}, and analogs of Ding's result for the CR-Yamabe problem on the Heisenberg group \cite{MaalaouiMartino2012, MaalaouiMartinoTralli2015}. Finally, we remark that in \cite{BahriChanillo2001, BahriXu2007} sign-changing solutions for critically nonlinear conformally invariant elliptic problems have been constructed using topological methods and in \cite{DelPino2011, DelPino2013}, sign-changing solutions for such problems have been constructed by arranging positive and negative bubbles (functions of the form \eqref{eq:scaled_translated_bubble}) whose locations and scales are carefully chosen.%

In the present work, the theme of extending Theorem \ref{oldtheorem:Ding} is continued, this time in the context of the conformally invariant Choquard equation \eqref{eq:choquard}. For particular values of $\mu$, the classification of non-negative solutions to problem \eqref{eq:choquard} was carried out in \cite{DuYang2019}. A corollary of their result is as follows: 
\begin{oldtheorem}
\label{oldtheorem:DuYang}
Let $n\geq 3$ and let $\mu\in (0, n)$. If one of the following holds 
\begin{enumerate}[label = {\bf(\alph*)}, ref = {\bf(\alph*)}, wide = 0pt]
	\item $n\in \{3, 4\}$ or \label{item:n_small}
	\item $n\geq 5$ and $\mu\in (0, 4]$  \label{item:n_large_mu_small}
\end{enumerate}
then there is a constant $a = a(n, \mu)>0$ such that for all positive solutions $u$ to \eqref{eq:choquard} there is $(x_0, \lambda)\in \bb R^n\times(0, \infty)$ for which $u = aU_{x_0, \lambda}$, where $U_{x_0, \lambda}$ is as in \eqref{eq:scaled_translated_bubble}.
\end{oldtheorem}
By using the fact that the functions $f_{x_0, \lambda} = U_{x_0, \lambda}^{2^*_\mu}$ with $U_{x_0, \lambda}$ as in \eqref{eq:scaled_translated_bubble} are extremal for the sharp Hardy-Littlewood-Sobolev inequality stated in the form of Theorem \ref{theorem:HLS} below or by using the technique of Lemma 5.6 in \cite{Li2004} one can directly verify that up to positive constant multiple and independently of the constraints imposed in items \ref{item:n_small} and \ref{item:n_large_mu_small} of Theorem \ref{oldtheorem:DuYang}, the functions \eqref{eq:scaled_translated_bubble} satisfy \eqref{eq:choquard} for all $\mu\in (0, n)$. However, to the authors' knowledge, the classification of positive solutions to problem \eqref{eq:choquard} for $n\geq 5$ and $\mu\in (4, n)$ is still open. The following theorem is the main result of this work.
\begin{theorem}
\label{theorem:main}
Let $n\geq 3$ and let $\mu\in (0, n)$. There is a sequence $(u_i)\subset D^{1, 2}(\bb R^n)$ of sign-changing solutions to \eqref{eq:choquard} for which $\|\Grad u_i\|_{L^2(\bb R^n)}\to \infty$ as $i\to\infty$. 
\end{theorem}

The proof of Theorem \ref{theorem:main} is based on lifting problem \eqref{eq:choquard} to the sphere and leveraging the symmetries of the sphere. Specifically, we will consider two isometric group actions on $\bb S^n$ and their induced actions on the Sobolev space $H^1(\bb S^n)$. For each of these actions we will produce an unbounded sequence of solutions to the lifted problem whose members are invariant under the induced action on $H^1(\bb S^n)$. The first action is inspired by Ding’s approach and it arises as the restriction of the standard action of $O(n + 1)$ to certain subgroups of $O(n + 1)$. When $n\geq 3$ and $\mu\in (0, n)$, this action always yields an unbounded sequence of solutions to the lifted problem. However, the deduction that these solutions are sign-changing relies on Theorem \ref{oldtheorem:DuYang} so we can only deduce the sign-changing behavior of these solutions if either $n\in \{3, 4\}$ or if both $n\geq 5$ and $\mu\in(0, 4]$ (corresponding to items \ref{item:n_small} and \ref{item:n_large_mu_small} of Theorem \ref{oldtheorem:DuYang} respectively). The second group action on $\bb S^n$ that we consider is inspired by the approach in \cite{BartschSchneiderWeth2004} and it arises as a non-standard action of certain subgroups of $O(n + 1)$ on $\bb S^n$. This method can be carried out for all $n\geq 3$ other than $n = 4$ but we do not know whether this technique can be carried out when $n = 4$. When this approach is applicable (e.g., when both $n\geq 3$ and $n\neq 4$), it produces an unbounded sequence of solutions to the lifted problem whose members are invariant under the action. Moreover, because of the definition of the group action, these solutions are always sign-changing independently of any further assumptions on $\mu$ or $n$. In fact, the second method also has the advantage that it allows one to draw some conclusions on the nodal sets of these solutions.

The paper is organized as follows. 
Section \ref{s:preliminaries} is devoted to the discussion of some preliminary notions including the functional space setting and the variational structure for problem \eqref{eq:choquard}. In Section \ref{s:lift_to_sphere} the variational form for problem \eqref{eq:choquard} is lifted via stereographic projection to a variational form defined on a Sobolev space over the standard sphere. In Section \ref{s:improved_embeddings} improved  embedding and compactness results for symmetric Sobolev spaces are used to show that the restriction of the lifted variational form to suitable subspaces of symmetric functions satisfies the Palais-Smale compactness condition. In Section \ref{s:existence} we use the Symmetric Mountain Pass Theorem of Ambrosetti and Rabinowitz and the Principle of Symmetric Criticality to prove Theorem \ref{theorem:main}. In Appendix \ref{s:appendix} an outline of a proof for the regularity of solutions to problem \eqref{eq:choquard} is provided.  
\section{Preliminaries}
\label{s:preliminaries}
We assume throughout the manuscript that $n\geq 3$. In this section we discuss the variational structure for problem \eqref{eq:choquard} and we reformulate this problem as an equivalent problem on $\bb S^n$. The space of smooth compactly supported functions on $\bb R^n$ is denoted by $C_c^\infty(\bb R^n)$ and the completion of $C_c^\infty(\bb R^n)$ relative to the norm $\|u\|_{D^{1, 2}(\bb R^n)} = \|\Grad u\|_{L^2(\bb R^n)}$ is denoted by $D^{1, 2}(\bb R^n)$.  It is well-known that $D^{1, 2}(\bb R^n)$ may alternatively characterized by 
\begin{equation*}
	D^{1, 2}(\bb R^n)
	= \{u\in L^{2^*}(\bb R^n): \|u\|_{D^{1, 2}(\bb R^n)}< \infty\}, 
\end{equation*}
where as in \eqref{eq:2_star}, $2^* = 2^*_0 = 2n/(n - 2)$. We will make extensive use of the Hardy-Littlewood-Sobolev inequality which we recall here. 
\begin{oldtheorem}[Hardy-Littlewood-Sobolev inequality]
\label{theorem:HLS}
Let $n\geq 1$ and let $\mu\in (0, n)$. For any $p, r\in(1, \infty)$ that satisfy $\frac 1p + \frac 1 r + \frac \mu n = 2$ there is a sharp constant $\mc H(n, \mu, p)> 0$ such that the inequality 
\begin{equation}
\label{eq:HLS_original}
	\abs{ \int_{\bb R^n}\int_{\bb R^n}\frac{f(y) h(x)}{|x - y|^\mu}\; \d y\; \d x} \leq \mc H(n, p, \mu)\|f\|_{L^p(\bb R^n)}\|h\|_{L^r(\bb R^n)} 
\end{equation}
holds for all $f\in L^p(\bb R^n)$ and all $h\in L^r(\bb R^n)$. 
\end{oldtheorem}
By Lebesgue duality, under the hypotheses of Theorem \ref{theorem:HLS} and with $\mc H(n, p, \mu)$ as in Theorem \ref{theorem:HLS}, one may equivalently express inequality \eqref{eq:HLS_original} as 
\begin{equation}
\label{eq:HLS_inequality}
	\|I_\mu f\|_{L^{r'}(\bb R^n)}\leq \mc H(n, \mu, p) \|f\|_{L^p(\bb R^n)}
\end{equation}
where $r' = r/(r - 1)$ is the Lebesgue conjugate of $r$. For brevity we will refer to each of inequalities \eqref{eq:HLS_original} and \eqref{eq:HLS_inequality} as the HLS inequality. 
In the diagonal case $p= r = 2n/(2n - \mu)$ Lieb showed in \cite{Lieb1983} that equality holds in \eqref{eq:HLS_inequality} if and only if $f = cU_{x_0, \lambda}^{2_\mu^*}$ for some $(c, x_0, \lambda)\in \bb R\times \bb R^n\times(0, \infty)$ where $U_{x_0, \lambda}$ as as in \eqref{eq:scaled_translated_bubble}. 
Equivalently, equality holds in \eqref{eq:HLS_original} if and only if there is $(x_0, \lambda)\in \bb R^n\times(0, \infty)$ for which each of $f$ and $h$ are constant multiples of $U_{x_0, \lambda}^{2^*_\mu}$. The following lemma shows that problem \eqref{eq:choquard} is critically nonlinear in the sense that the right-hand side of the equation lies in the Lebesgue dual of $L^{2^*}(\bb R^n)$. 
\begin{lemma}
\label{lemma:rhs_is_critical}
Let $n\geq 3$ and let $\mu\in (0, n)$. There is a constant $C = C(n, \mu)> 0$ such that 
\begin{equation*}
	\|I_\mu[|u|^{2^*_\mu}]|u|^{2^*_\mu - 2}u\|_{L^{2n/(n + 2)}(\bb R^n)}
	\leq C\|u\|_{L^{2^*}(\bb R^n)}^{2\cdot 2^*_\mu - 1}
	\qquad \text{ for all } u\in L^{2^*}(\bb R^n).  
\end{equation*}
\end{lemma}
\begin{proof}
Applying inequality \eqref{eq:HLS_inequality} with $p = 2^*/2^*_\mu$ and $r' = 2n/\mu$ gives a constant $C = C(n, \mu)>0$ such that
\begin{equation*}
	\|I_\mu[|u|^{2^*_\mu}]\|_{L^{2n/\mu}(\bb R^n)}
	\leq C\||u|^{2^*_\mu}\|_{L^{2^*/2^*_\mu}(\bb R^n)}
	= C\|u\|_{L^{2^*}(\bb R^n)}^{2^*_\mu} 
\end{equation*}
for all $u\in L^{2^*}(\bb R^n)$. For any such $u$, using this estimate together with H\"older's inequality gives
\begin{equation*}
\begin{split}
	\|I_\mu[|u|^{2^*_\mu}]|u|^{2^*_\mu - 2}u\|_{L^{2n/(n + 2)}(\bb R^n)}
	& \leq \|I_\mu[|u|^{2^*_\mu}]\|_{L^{2n/\mu}(\bb R^n)}\|u\|_{L^{2^*}(\bb R^n)}^{2^*_\mu - 1}\\
	 & \leq C(n, \mu)\|u\|_{L^{2^*}(\bb R^n)}^{2\cdot2^*_\mu - 1}.
\end{split}
\end{equation*}
\end{proof}
In view of Lemma \ref{lemma:rhs_is_critical} and the Sobolev embedding $D^{1, 2}(\bb R^n)\hookrightarrow L^{2^*}(\bb R^n)$, problem \eqref{eq:choquard} may be understood in the weak sense. 
\begin{defn}
A \emph{weak solution} to problem \eqref{eq:choquard} is any function $u\in D^{1, 2}(\bb R^n)$ for which 
\begin{equation*}
	\int_{\bb R^n}\Grad u\cdot \Grad \varphi = \int_{\bb R^n}I_\mu[|u|^{2^*_\mu}]|u|^{2^*_\mu - 2} u\varphi
	\qquad \text{ for all }\varphi \in D^{1, 2}(\bb R^n). 
\end{equation*}
\end{defn}
The next theorem provides the regularity properties of weak solutions to problem \eqref{eq:choquard}. In the statement of the theorem the containment $u\in C^{2, \alpha}(\bb R^n)$ is understood to mean that $\|u\|_{C^{2, \alpha}(K)}< \infty$ for all compact subsets $K\subset \bb R^n$. 
\begin{theorem}
\label{theorem:regularity}
There exists $\alpha\in (0, 1)$ such that every weak solution $u$ to problem \eqref{eq:choquard} satisfies $u\in C^{2, \alpha}(\bb R^n)$. 
\end{theorem}
Although the proof of Theorem \ref{theorem:regularity} can be obtained by standard (yet lengthy) techniques, the authors were unable to find a proof in the literature. For convenience a proof is outlined in Appendix \ref{s:appendix}.

By performing standard computations involving the Sobolev inequality and the HLS inequality, one can verify that the functional $F:D^{1, 2}(\bb R^n)\to \bb R$ given by 
\begin{equation}
\label{eq:F_functional}
	F(u) = \frac 12\int_{\bb R^n}|\Grad u|^2 - \frac{1}{2\cdot2^*_\mu}\int_{\bb R^n}I_\mu[|u|^{2^*_\mu}]|u|^{2^*_\mu}
\end{equation}
is of class $C^1$. Moreover, the weak solutions to problem \eqref{eq:choquard} are precisely the critical points of $F$. Thus, our approach in proving Theorem \ref{theorem:main} is to establish the existence of a sequence of sign-changing critical points $(u_i)_{i = 1}^\infty\subset D^{1, 2}(\bb R^n)$ of $F$ for which $\|u_i\|_{D^{1, 2}(\bb R^n)}\to\infty$. The primary obstacle in doing so is the failure of compactness in the embedding $D^{1, 2}(\bb R^n)\hookrightarrow L^{2^*}(\bb R^n)$. 
\section{An Equivalent Problem on $\bb S^n$}
\label{s:lift_to_sphere}
In this section we use stereographic projection to lift the functional $F$ defined in equation \eqref{eq:F_functional} to an equivalent functional defined on a Sobolev space over the standard sphere $\bb S^n$. Let $S = (0, \ldots, 0, -1)\in \bb R^{n + 1}$ denote the south pole of $\bb S^n$ and let $\pi:\bb S^n\setminus\{S\}\to \bb R^n$ be stereographic projection through the south pole to the equatorial plane of $\bb R^n$ so that  
\begin{equation}
\label{eq:stereographic_projection}
	\pi(\xi) = \left(\frac{\xi_1}{1 + \xi_{n + 1}}, \ldots, \frac{\xi_n}{1 + \xi_{n + 1}}\right)
	\qquad \xi \in \bb S^n\setminus\{S\} 
\end{equation}
and 
\begin{equation*}
	\pi^{-1}(x) = \left(\frac{2x}{1 + |x|^2}, \frac{1 - |x|^2}{1 + |x|^2}\right)
	\qquad x\in \bb R^n. 
\end{equation*}
Stereographic projection induces the ``projection'' map $P:C_c^\infty(\bb S^n\setminus\{S\})\to C_c^\infty(\bb R^n)$ at the level of functional spaces defined by 
\begin{equation}
\label{eq:P_projection}
	Pv(x) = \left(\frac 2{1 + |x|^2}\right)^{\frac{n - 2}{2}} v\circ \pi^{-1}(x) \qquad \text{ for }x\in \bb R^n.
\end{equation}
This map evidently both linear and invertible and the inverse of $P$ is given by
\begin{equation*}
	P^{-1} u(\xi) = \left(\frac{1}{1 + \xi_{n + 1}}\right)^{\frac{n - 2}{2}} u\circ \pi(\xi)
	\qquad \text{ for } \xi \in \bb S^n. 
\end{equation*}
One can routinely verify that $P$ preserves $L^{2^*}$-norms in the sense that $\|Pv\|_{L^{2^*}(\bb R^n)} = \|v\|_{L^{2^*}(\bb S^n)}$ whenever $v\in C_c^\infty(\bb S^n\setminus\{S\})$. Moreover, letting $H^1(\bb S^n)$ denote the completion of $C^\infty(\bb S^n)$ relative to the norm 
\begin{equation}
\label{eq:correct_H1_norm}
	\|v\|_{H^1(\bb S^n)}^2
	= \int_{\bb S^n}\left(|\Grad v|^2 + \frac{n(n - 2)}4 v^2\right)\; \d V,  
\end{equation} 
where $\d V$, $\Grad$ and $|\cdot|^2$ are all defined relative to the round metric on $\bb S^n$, one can also routinely verify that $\|Pv\|_{D^{1, 2}(\bb R^n)} = \|v\|_{H^1(\bb S^n)}$ for all $v\in C_c^\infty(\bb S^n\setminus\{S\})$. The proof of the following lemma is standard and therefore omitted. 
\begin{lemma}
\label{lemma:extend_P}
The map $P:C_c^\infty(\bb S^n\setminus\{S\})\to C_c^\infty(\bb R^n)$ defined in \eqref{eq:P_projection} extends to a norm-preserving bijective linear map $L^{2^*}(\bb S^n)\to L^{2^*}(\bb R^n)$ whose restriction $P|_{H^1(\bb S^n)}: H^1(\bb S^n)\to D^{1, 2}(\bb R^n)$ is also norm-preserving. 
\end{lemma}
For ease of notation we continue to denote the extension of $P$ (as described in Lemma \ref{lemma:extend_P}) by $P$. With the definition and elementary properties of $P$ in hand, we proceed to lift problem \eqref{eq:choquard} to an equivalent problem on $\bb S^n$. Specifically, since problem \eqref{eq:choquard} is the Euler-Lagrange equation for the functional $F$ given in \eqref{eq:F_functional}, we will show that the critical points of the functional $E:H^1(\bb S^n)\to \bb R$ given by 
\begin{equation}
\label{eq:E_functional}
	E = F\circ P
\end{equation} 
are in one-to-one correspondence with those of $F$. 
The main result of this section is the following proposition with clarifies the sense in which problem \eqref{eq:choquard} and it's lift to the sphere are equivalent. 
\begin{prop}
\label{prop:critical_point_correspondence}
A function $u\in D^{1, 2}(\bb R^n)$ is a critical point of the functional $F:D^{1, 2}(\bb R^n)\to \bb R$ given in \eqref{eq:F_functional} if and only if $P^{-1}u$ is a critical point of the functional $E$ defined in \eqref{eq:E_functional}. In particular, the critical points of $F$ are in one-to-one correspondence with those of $E$.
\end{prop}
The proof of Proposition \ref{prop:critical_point_correspondence} will be established with the aid of a series of lemmata. We begin by deriving an explicit expression for $E(v)$ whenever $v\in H^1(\bb S^n)$, see Lemmata \ref{lemma:derivative_Imu_stuff} and \ref{lemma:E_explicit} that follow. 
\begin{lemma}
\label{lemma:derivative_Imu_stuff}
Let $n\geq 3$ and let $\mu\in (0, n)$. For every $u, \psi\in L^{2^*}(\bb R^n)$ the equality 
\begin{equation*}
	\int_{\bb R^n}I_\mu[|u|^{2^*_\mu}]|u|^{2^*_\mu - 2}u \psi \; \d x
	= \int_{\bb S^n}J_\mu[|v|^{2^*_\mu}]|v|^{2^*_\mu - 2}v\varphi \; \d V
\end{equation*}
holds, where $v, \varphi\in L^{2^*}(\bb S^n)$ are defined by $v = P^{-1}u$ and $\varphi = P^{-1}\psi$ and where $J_\mu$ is defined by 
\begin{equation}
\label{eq:J_mu}
	J_\mu f(\xi) = \int_{\bb S^n}\frac{f(\zeta)}{|\xi- \zeta|^\mu}\; \d V(\zeta)
\end{equation}
with $|\xi- \zeta|$ denoting the chordal distance in $\bb R^{n + 1}$ between $\xi, \zeta\in \bb S^n$. 
\end{lemma}
\begin{proof}
With $\pi:\bb S^n\setminus\{S\}\to \bb R^n$ as in \eqref{eq:stereographic_projection}, for any $\xi, \zeta\in \bb S^n$, the distance in $\bb R^n$ between $\pi(\xi)$ and $\pi(\zeta)$ is related to the chordal distance $|\xi - \zeta|$ in $\bb R^{n + 1}$ via 
\begin{equation*}
	|\pi(\xi) - \pi(\zeta)|
	= \frac{|\xi - \zeta|}{\sqrt{1 + \xi_{n + 1}}\sqrt{1 + \zeta_{n + 1}}}. 
\end{equation*}
Therefore, for any $\xi\in \bb S^n$, using the change of variable $y = \pi(\zeta)$ we have $\d y = (1 + \zeta_{n + 1})^{-n}\; \d V(\zeta)$ and 
\begin{equation*}
\begin{split}
	I_\mu[|u|^{2^*_\mu}](\pi(\xi))
	&= \int_{\bb R^n}\frac{|u(y)|^{2^*_\mu}}{|\pi(\xi) - y|^\mu}\; \d y\\
	& = \int_{\bb S^n}\frac{|u(\pi(\zeta))|^{2^*_\mu}}{|\pi(\xi) - \pi(\zeta)|^\mu}\; \left(\frac 1{1 + \zeta_{n+ 1}}\right)^n\; \d V(\zeta)\\ 
	& = (1 + \xi_{n + 1})^{\mu/2}\int_{\bb S^n}\frac{|P^{-1}u(\zeta)|^{2^*_\mu}}{|\xi - \zeta|^\mu}\; \d V(\zeta)\\
	& = (1+ \xi_{n + 1})^{\mu/2}J_\mu[|v|^{2^*_\mu}](\xi). 
\end{split}
\end{equation*}
Using this equality together with the change of variable $x = \pi(\xi)$ 
we have 
\begin{equation*}
\begin{split}
	\int_{\bb R^n}& I_\mu[|u|^{2^*_\mu}](x)|u(x)|^{2^*_\mu - 2}u(x) \psi(x) \; \d x\\
	& = \int_{\bb S^n}I_\mu[|u|^{2^*_\mu}](\pi(\xi))|u(\pi(\xi))|^{2^*_\mu - 2}u(\pi(\xi)) \psi(\pi(\xi)) \; \d V(\xi)\\
	& = \int_{\bb S^n}\left(\frac 1{1 + \xi_{n + 1}}\right)^{\mu/2}I_\mu[|u|^{2^*_\mu}](\pi(\xi))|v(\xi)|^{2^*_\mu - 2}v(\xi) \varphi(\xi) \; \d V(\xi)\\
	& = \int_{\bb S^n}J_\mu[|v|^{2^*_\mu}](\xi)|v(\xi)|^{2^*_\mu - 2}v(\xi)\varphi(\xi)\; \d V(\xi). 
\end{split}
\end{equation*}
\end{proof}
The following lemma provides an explicit expression for $E(v)$ whenever $v\in H^1(\bb S^n)$. 
\begin{lemma}
\label{lemma:E_explicit}
Let $F$ be as in \eqref{eq:F_functional} and let $E$ be as in \eqref{eq:E_functional}. For any $v\in H^1(\bb S^n)$ the equality 
\begin{equation}
\label{eq:E_explicit_expression}
	E(v) = \frac 12\int_{\bb S^n}\left(|\Grad v|^2 + \frac{n(n- 2)}4v^2\right)\; \d V- \frac 1{2\cdot 2^*_\mu}\int_{\bb S^n}J_\mu[|v|^{2^*_\mu}]|v|^{2^*_\mu}\; \d V 
\end{equation}
holds, where $J_\mu$ is as in equation \eqref{eq:J_mu}. 

\end{lemma}
\begin{proof}
The asserted equality follows immediately from the fact that $P:H^1(\bb S^n)\to D^{1, 2}(\bb R^n)$ is norm-preserving (see Lemma \ref{lemma:extend_P}) and from Lemma \ref{lemma:derivative_Imu_stuff} applied with $u = \psi = Pv$. 
\end{proof}
\begin{proof}[Proof of Proposition \ref{prop:critical_point_correspondence}]
In view of the bijectivity of $P:H^1(\bb S^n)\to D^{1, 2}(\bb R^n)$, it suffices to show that 
\begin{equation*}
	\lb F'(u), \psi\rb = \lb E'(v), \varphi\rb
\end{equation*}
for all $u, \psi\in D^{1, 2}(\bb R^n)$, where $v = P^{-1}u$ and $\varphi = P^{-1}\psi$. Here we use $\lb \cdot, \cdot \rb$ to denote both the dual pairing of $[D^{1, 2}(\bb R^n)]'$ with $D^{1, 2}(\bb R^n)$ and the dual pairing of $[H^1(\bb S^n)]'$ with $H^1(\bb S^n)$. In each usage of this notation, the context makes it clear which dual pairing is intended. Now fixing any $u,\psi\in D^{1, 2}(\bb R^n)$, in view of the polarization identity for real inner products and the fact that $P:H^1(\bb S^n)\to D^{1, 2}(\bb R^n)$ is both linear and norm-preserving we have
\begin{equation*}
\begin{split}
	2\int_{\bb R^n}\Grad u\cdot \Grad\psi
	& = \|u\|_{D^{1, 2}(\bb R^n)}^2 + \|\psi\|_{D^{1, 2}(\bb R^n)}^2 - \|u- \psi\|_{D^{1, 2}(\bb R^n)}^2\\
	& = \|v\|_{H^1(\bb S^n)}^2 + \|\varphi\|_{H^1(\bb S^n)}^2 - \|v- \varphi\|_{H^1(\bb S^n)}^2\\
	& = 2\int_{\bb S^n}\left(g_{\bb S^n}(\Grad v, \Grad \varphi) + \frac{n(n - 2)}4 v\varphi\right)\; \d V. 
\end{split}	
\end{equation*}
Combining this equality with Lemma \ref{lemma:derivative_Imu_stuff} we have
\begin{equation*}
\begin{split}
	\lb F'(u), \psi\rb
	& = \int_{\bb R^n}\Grad u\cdot \Grad\psi - \int_{\bb R^n}I_\mu[|u|^{2^*_\mu}]|u|^{2^*_\mu - 2}u \psi\\
	& = \int_{\bb S^n}\left(g_{\bb S^n}(\Grad v,\Grad \varphi) + \frac{n(n - 2)}4 v\varphi\right)\; \d V - \int_{\bb S^n}J_\mu[|v|^{2^*_\mu}]|v|^{2^*_\mu - 2}\varphi\; \d V\\
	& = \lb E'(v), \varphi\rb. 
\end{split}
\end{equation*}
\end{proof} 
\section{Improved Embeddings and Compactness}
\label{s:improved_embeddings}

In this section we discuss two isometric group actions on $H^1(\bb S^n)$ and corresponding improved Sobolev embeddings and compactness results for the subspaces of $H^1(\bb S^n)$ that are invariant under these actions. If $G$ is a subgroup of isometries of a Riemannian manifold $(M, g)$ and if $\xi\in M$ we denote the orbit of $\xi$ under the action of $G$ by $G\xi = \{\gamma(\xi): \gamma\in G\}$ and we denote by $H_G^1(M)$ the subspace of $H^1(M)$ consisting of the functions in $H^1(M)$ that are invariant under the action of $G$. The following theorem gives the case of the improved Sobolev embedding in the presence of symmetries that will be needed in the sequel. One may consult Theorem 9.1 of \cite{Hebey2000} for a more general statement of the theorem and for the proof of the theorem.
\begin{oldtheorem}
\label{oldtheorem:improved_embedding}
Let $(M, g)$ be a smooth compact Riemannian $n$-manifold and let $G$ be a compact subgroup of isometries of $M$. If $\card G\xi = +\infty$ for all $\xi\in M$ then for $d$ defined by 
\begin{equation}
\label{eq:minimum_orbit_dimension}
	d = \min\{\dim G\xi: \xi\in M\}
\end{equation} 
we have $d\geq 1$ and 
\begin{enumerate}[label = {\bf(\alph*)}, wide = 0pt]
	\item If $n - d\leq 2$ then for any $p\geq 1$ $H_G^1(M)\subset L^p(M)$ and the embedding is both continuous and compact. 
	\item If $n - d> 2$ then for any $p\in [1, \frac{2(n - d)}{n - d - 2}]$ the embedding $H_G^1(M)\subset L^p(M)$ is continuous. This embedding is compact provided $p\in [1, \frac{2(n - d)}{n - d - 2})$. 
\end{enumerate}
In particular, there exists $p_0> 2^*$ such that, for any $p\in [1, p_0]$ the embedding $H_G^1(M)\subset L^p(M)$ is continuous and compact. 
\end{oldtheorem}
We now describe the two group actions on $\bb S^n$ and their associated actions on $H^1(\bb S^n)$ that will be used to construct solutions to the lift of problem \eqref{eq:choquard} to the sphere. The first action is the standard action on $H^1(\bb S^n)$ induced by the action of $O(n + 1)$ on $\bb S^n$. Namely, the action of $h\in O(n + 1)$ on $v\in H^1(\bb S^n)$ is 
\begin{equation}
\label{eq:standard_action}
	h\cdot v = v\circ h^{-1}. 
\end{equation}
To describe the second action, we make the following definition. In the statement of the definition $NG$ denotes the normalizer in $O(n + 1)$ of a subgroup $G\subset O(n + 1)$. 
\begin{defn}
\label{defn:property_P}
A subgroup $G\subset O(n + 1)$ is said to satisfy property $\mathcal P$ if there exists $\tau\in NG\setminus G$ for which both $\tau^2\in G$ and 
\begin{equation}
\label{eq:tau_orbit}
	\tau\xi\not\in G\xi \qquad \text{ for some }\xi\in \bb S^n. 
\end{equation}
\end{defn}
The standard example of a subgroup $G\subset O(n+ 1)$ and an element $\tau\in NG\setminus G$ for which property $\mathcal P$ holds is as follows.
\begin{example}
\label{example:standard}
For $m\geq 2$, let 
\begin{equation}
\label{eq:G_form}
	G = O(n_1)\times O(n_2)\times \ldots\times O(n_m)
\end{equation}
for some positive integers $n_1, \ldots, n_m$ satisfying both 
\begin{equation}
\label{eq:nj_sum}
	\sum_{j = 1}^m n_j = n + 1
\end{equation} 
and $n_1 = n_2$. Define $\tau$ by its matrix representation relative to the standard basis of $\bb R^{n + 1}$:
\begin{equation}
\label{eq:standard_tau}
	\tau = \begin{bmatrix}
	0 & I_{n_1} & 0\\
	I_{n_1}& 0 & 0\\
	0 & 0 & I_{n + 1-2n_1}
	\end{bmatrix}. 
\end{equation}
Evidently we have both $\tau\in NG\setminus G$ and $\tau^2\in G$. Moreover, if $\xi = (\xi_1, \xi_2, \xi_3)\in \bb S^n\subset \bb R^{n_1}\times \bb R^{n_1}\times \bb R^{n + 1 - 2n_1}$ and if $|\xi_1|\neq|\xi_2|$ then $\tau\xi\not\in G\xi$. 
\end{example}
The following lemma indicates the relevance of the conditions $\tau\in NG\setminus G$ and $\tau^2\in G$ in property $\mathcal P$. In the statement of the lemma $L(H^1(\bb S^n), H^1(\bb S^n))$ denotes the space of bounded linear maps on $H^1(\bb S^n)$. Since the proof is elementary, we omit the details. 
\begin{lemma}
\label{lemma:admits_gamma}
If $G\subset O(n+ 1)$ is a subgroup for which there exists $\tau\in NG\setminus G$ that satisfies $\tau^2\in G$ then for any such $\tau$
\begin{equation}
\label{eq:Gamma_Gtau}
	\Gamma(G, \tau):= G\cup G\tau
\end{equation}
is a subgroup of $O(n + 1)$ in which $G$ is normal. Moreover, $\Gamma(G, \tau)$ admits the semidirect product decomposition 
\begin{equation*}
	\Gamma(G, \tau) 
	= G\rtimes \lb\tau\rb
	\cong G\rtimes O(1). 
\end{equation*}
Finally, the action of $\Gamma(G, \tau)$ on $\bb S^n$ induces an action $\beta:\Gamma(G, \tau)\to L(H^1(\bb S^n), H^1(\bb S^n))$ such that for every $\gamma\in \Gamma(G, \tau)$ the map $\beta(\gamma)\in L(H^1(\bb S^n), H^1(\bb S^n))$ is norm-preserving. This action is completely determined by the equalities 
\begin{equation}
\label{eq:non_standard_action}
\begin{cases}
	\beta(h)v = v\circ h^{-1} & \text{ for }h\in G\\
	\beta(\tau)v = -v\circ \tau^{-1}  
\end{cases}
\end{equation}
for  $v\in H^1(\bb S^n)$.
\end{lemma}
For a subgroup $G\subset O(n + 1)$ let us denote by 
\begin{equation}
\label{eq:H1G}
	H^1_G(\bb S^n) = \{v\in H^1(\bb S^n): v \circ h^{-1}= v\text{ for all }h\in G\}
\end{equation}
the subspace of $H^1(\bb S^n)$ that is invariant under the (standard) action of $G$. If $G$ satisfies property $\mathcal P$ with corresponding element $\tau\in NG\setminus G$, let us denote by 
\begin{equation*}
	H^1_{\Gamma(G, \tau)}(\bb S^n) 
	= \{v\in H^1_G(\bb S^n): -v\circ\tau^{-1} = v\} 
\end{equation*}
the subspace of $H^1(\bb S^n)$ consisting of those functions that are invariant relative to the (non-standard) action $\beta$ as defined in \eqref{eq:non_standard_action}. Observe that if one omits condition \eqref{eq:tau_orbit} in property $\mathcal P$ then $H_{\Gamma(G, \tau)}^1(\bb S^n) = \{0\}$. Indeed, in the absence of this condition, for every $\xi\in \bb S^n$ there is $h = h_\xi\in G$ such that $\tau\xi = h\xi$. If $v\in H_{\Gamma(G, \tau)}^1(\bb S^n)$ and $\xi\in \bb S^n$, then with $h = h_\xi$ as such we have 
\begin{equation*}
	-v(\xi)
	= v\circ \tau(\xi)
	= v\circ h(\xi)
	= v(\xi)
\end{equation*}
and thus $v\equiv 0$. 

We are particularly interested in subgroups $G\subset O(n + 1)$ of the form \eqref{eq:G_form} that satisfy \eqref{eq:nj_sum} and the additional assumption
\begin{equation}
\label{eq:nj_largeness}
	\min_j\{n_j\} \geq 2. 
\end{equation} 
For any such $G$, we have $\card G \xi= +\infty$ for all $\xi\in \bb S^n$. If, in addition, $G$ satisfies property $\mathcal P$ with corresponding element $\tau\in NG\setminus G$, then the containment $G\subset\Gamma (G, \tau)$ guarantees that $\card \Gamma\xi\geq \card G\xi =+\infty$ for all $\xi\in \bb S^n$. In particular we obtain the following corollary of Theorem \ref{oldtheorem:improved_embedding}. 
\begin{coro}
\label{coro:improved_embedding}
Let $n\geq 3$, let $G$ be as in \eqref{eq:G_form} for some integer $m\geq 2$ and some integers $n_1, \ldots, n_m$ satisfying both \eqref{eq:nj_sum} and \eqref{eq:nj_largeness}. There exists $p>2^*$ and a constant $C(p)> 0$ such that 
\begin{equation}
\label{eq:improved_embedding}
	\|v\|_{L^p(\bb S^n)}\leq C\|v\|_{H^1(\bb S^n)}
	\qquad \text{ for all }v\in H_G^1(\bb S^n) 
\end{equation}
and for which the embedding $H_G^1(\bb S^n)\hookrightarrow L^p(\bb S^n)$ is compact. If, in addition to the above hypotheses $G$ satisfies property $\mathcal P$ with corresponding element $\tau\in NG\setminus G$, then (for the same $p$) \eqref{eq:improved_embedding} holds for all $v\in H_{\Gamma(G, \tau)}^1(\bb S^n)$ and the embedding $H_{\Gamma(G, \tau)}^1(\bb S^n)\hookrightarrow L^p(\bb S^n)$ is compact. 
\end{coro}
\begin{remark}
If $n = 4$, then there is no partition of $n + 1 = 5$ consisting of at least two elements for which two elements are equal and for which no element is less than $2$. In this case the construction of Example \ref{example:standard} does not yield a subgroup $G\subset O(5)$ for which property $\mathcal P$ holds. 
\end{remark}
\begin{lemma}
\label{lemma:E_invariant}
The functional $E$ defined in \eqref{eq:E_functional} is invariant under the standard action of $O(n + 1)$ defined in \eqref{eq:standard_action} in the sense that for any $v\in H^1(\bb S^n)$ and any $h\in O(n + 1)$, the equality 
\begin{equation}
\label{eq:E_invariant}
	E(hv) =E(v)
\end{equation} 
holds. If $G$ satisfies property $\mathcal P$ with corresponding element $\tau\in NG\setminus G$ then $E$ is invariant under the action of $\Gamma(G, \tau)$ on $H^1(\bb S^n)$ as defined in \eqref{eq:non_standard_action}. 
\end{lemma}
\begin{proof}
In view of the definition of the action of $\Gamma(G, \tau)$ on $H^1(\bb S^n)$ and the fact that $E$ is even, it suffices to show that $E$ is invariant under the standard action of $O(n + 1)$ on $H^1(\bb S^n)$. To do so we separately verify the invariance under the action of $O(n + 1)$ for each term in the explicit form of $E$ given in \eqref{eq:E_explicit_expression}. Accordingly, fix $h\in O(n + 1)$. To handle the first term in \eqref{eq:E_explicit_expression}, we use the equality $|\Grad(v\circ h^{-1})|^2 = |\Grad v|^2\circ h^{-1}$, the change of variable $\zeta\mapsto h\zeta$, and the fact that $\d V$ is invariant under the action of $h$ as follows: 
\begin{equation}
\label{eq:term1_invariance}
\begin{split}
	\int_{\bb S^n}&  \left(|\Grad (h v)|^2 + \frac{n(n - 2)}4(hv)^2\right)\; \d V\\
	= &\; \int_{\bb S^n}\left(|\Grad v|^2\circ h^{-1} + \frac{n(n - 2)}4 v^2\circ h^{-1}\right)\; \d V\\
	= &\;  \int_{\bb S^n}\left(|\Grad v|^2 + \frac{n(n - 2)}4 v^2\right)\; \d V. 
\end{split}
\end{equation}
To handle the second term in \eqref{eq:E_explicit_expression}, we first use the change of variable $\zeta \mapsto h\zeta$ and the fact that the volume form on $\bb S^n$ is invariant under the action of $h$ to find that for any $w\in L^{2^*/2^*_\mu}(\bb S^n)$,
\begin{equation*}
\begin{split}
	J_\mu[h w](h\xi)
	& = \int_{\bb S^n}\frac{w\circ h^{-1}(\zeta)}{|h\xi - \zeta|^\mu}\; \d V(\zeta)\\
	& = \int_{\bb S^n} \frac{w(\zeta)}{|h\xi - h\zeta|^\mu}\; \d V(\zeta)\\
	& = J_\mu[w](\xi),  
\end{split}
\end{equation*}
the final equality holding as the action of $O(n+ 1)$ on $\bb S^n$ preserves the chordal distance between points of $\bb S^n$. Using this equality and using once more the change of variable $\zeta \mapsto h\zeta$ we have
\begin{equation}
\label{eq:term2_invariance}
\begin{split}
	\int_{\bb S^n}& J_\mu[|h v|^{2^*_\mu}](\zeta)|h v(\zeta)|^{2^*_\mu}\; \d V(\zeta)\\
	& = \int_{\bb S^n} J_\mu[|hv|^{2^*_\mu}](h\zeta)|h v(h\zeta)|^{2^*_\mu}\; \d V(\zeta)\\
	& = \int_{\bb S^n}J_\mu[|v|^{2^*_\mu}](\zeta)|v(\zeta)|^{2^*_\mu}\; \d V(\zeta). 
\end{split}
\end{equation}
Equation \eqref{eq:E_invariant} follows immediately from equalities \eqref{eq:term1_invariance} and \eqref{eq:term2_invariance} and from the explicit expression for $E$ given in \eqref{eq:E_explicit_expression}. 
\end{proof}
For $i\in \bb N\cup\{0\}$ we denote the space of spherical harmonics of degree $i$ by $\bb Y_i$ and for $G\subset O(n+ 1)$ we define 
\begin{equation}
\label{eq:symmeric_spherical_harmonics}
	\bb Y_{G, i}
	= \{Y\in \bb Y_i: h Y= Y\text{ for all }h\in G\}  
\end{equation}
to be the space of $G$-invariant spherical harmonics of degree $i$. If in addition $G$ satisfies property $\mathcal P$ with corresponding element $\tau\in NG\setminus G$ then with $\Gamma = \Gamma(G, \tau)$ as in \eqref{eq:Gamma_Gtau} we define
\begin{equation}
\label{eq:Gamma_symmeric_spherical_harmonics}
	\bb Y_{\Gamma, i}
	= \{Y\in \bb Y_{G, i}: Y\circ \tau^{-1}= -Y\}
\end{equation}
to be the space of $\Gamma$-invariant spherical harmonics of degree $i$.
\begin{lemma}
For any subgroup $G\subset O(n + 1)$, the symmetric Sobolev space $H_G^1(\bb S^n)$ defined in \eqref{eq:H1G} admits the following direct sum decomposition 
\begin{equation}
\label{eq:H1G_decomposition}
	H_G^1(\bb S^n) = \overline{\bigoplus_{i = 0}^\infty \bb Y_{G, i}}
\end{equation}
in the sense that for all $v\in H_G^1(\bb S^n)$, there is a sequence $(Y_i)_{i = 1}^\infty$ for which both $Y_i\in \bb Y_{G,i}$ for all $i$ and $\|v - \sum_{i = 0}^NY_i\|_{H^1(\bb S^n)}\to 0$ as $N\to\infty$. If, in addition, $G$ satisfies property $\mathcal P$ with corresponding element $\tau\in NG\setminus G$ then with $\Gamma = \Gamma(G, \tau)$ as in \eqref{eq:Gamma_Gtau} there holds
\begin{equation}
\label{eq:H1_Gamma_decomposition}
	H_{\Gamma}^1(\bb S^n) = \overline{\bigoplus_{i = 0}^\infty \bb Y_{\Gamma, i}}
\end{equation}
\end{lemma}
\begin{proof}
To show that decomposition \eqref{eq:H1G_decomposition} holds, let $G$ be any subgroup of $O(n + 1)$, let $v\in H_G^1(\bb S^n)$ and choose a sequence $(Y_i)_{i = 1}^\infty$ for which $Y_i\in \bb Y_i$ and $v = \sum_{i = 0}^\infty Y_i$. We need to show that each summand on the right-hand side of this equality is invariant under the action of $G$ in the sense that $Y_i = h Y_i$ for all $h\in G$ and all $i\in \{0, 1, 2,  \ldots\}$. Since $v$ is invariant under the action of $G$, for all $h\in G$ we have
\begin{equation*}
	\sum_{i = 0}^\infty Y_i
	= v
	= v\circ h^{-1}
	= \sum_{i = 0}^\infty Y_i\circ h^{-1}
	= \sum_{i = 0}^\infty hY_i. 
\end{equation*}
Considering fixed but arbitrary $h\in G$ and defining $Z_i = Y_i - h Y_i$, the above string of equalities gives
\begin{equation*}
	0 = \sum_{i = 0}^\infty(Y_i - h Y_i)
	= \sum_{i = 0}^\infty Z_i. 
\end{equation*}
For every $i\in\{0, 1, 2, \ldots\}$ the space $\bb Y_i$ is invariant under the action of $O(n + 1)$ and hence it is also invariant under the action of $G$. Thus, $h Y_i\in \bb Y_i$ for every $h\in G$ and every $Y_i\in \bb Y_i$ (but based only on the $G$-invariance, we do not have $h Y_i = Y_i$). In particular, $Z_i\in \bb Y_i$ for all $i$. By the orthogonality of the spherical harmonics we have $\int_{\bb S^n}Z_iZ_j\; \d V = 0$ whenever $i\neq j$. Therefore, as $N\to \infty$ we have 
\begin{equation*}
\begin{split}
	\circ(1)
	& = \|\sum_{i = 0}^NZ_i\|_{H^1(\bb S^n)}\\
	& \geq \frac{n(n - 2)}4\int_{\bb S^n}|\sum_{i = 0}^N Z_i|^2\; \d V\\
	& = \frac{n(n - 2)}4\sum_{i = 0}^N\int_{\bb S^n}Z_i^2\; \d V,
\end{split}
\end{equation*}
from which we deduce that $Z_i\equiv 0$ for all $i$. The proof that decomposition \eqref{eq:H1_Gamma_decomposition} holds whenever $G$ satisfies property $\mathcal P$ with corresponding element $\tau\in NG\setminus G$ is essentially the same as the proof of decomposition \eqref{eq:H1G_decomposition}. 
\end{proof}
\begin{lemma}
\label{lemma:weakly_closed}
Let $(M, g)$ be a closed compact Riemannian manifold and let $G$ be a subgroup of isometries of $M$. The space $H_G^1(M)$ of $G$-invariant functions is a weakly closed subspace of $H^1(M)$ in the sense that if $(v_i)_{i = 1}^\infty\subset H_G^1(M)$ and if $v\in H^1(M)$ with $v_i\rightharpoonup v$ weakly in $H^1(M)$, then $v\in H_G^1(M)$. 
\end{lemma}
\begin{proof}
Let $(v_i)_{i = 1}^\infty\subset H_G^1(M)$ and $v\in H^1(M)$ with $v_i\weakconv v$ weakly in $H^1(M)$ and let $h:M\to M$ be an isometry. For any $\varphi\in L^2(M)$ we identify both $\varphi$ and $h^{-1}\varphi$ with elements of $[H^1(M)]'$ and use both the change of variable $\zeta\mapsto h\zeta$ and the $G$-invariance of $v_i$ to obtain
\begin{equation*}
\begin{split}
	\int_M(h v) \; \varphi\; \d V
	& = \int_M v\; (h^{-1}\varphi)\; \d V\\
	& = \lim_i\int_M v_i\; (h^{-1}\varphi)\; \d V\\
	& = \lim_i\int_M (h v_i)\; \varphi\; \d V\\
	& = \lim_i\int_M v_i\varphi\; \d V\\
	& = \int_M v\varphi\; \d V. 
\end{split}
\end{equation*}
Since $\varphi\in L^2(M)$ is arbitrary we conclude that $\|h v - v\|_{L^2(M)} = 0$. It remains to show that 
\begin{equation}
\label{eq:gradient_symmetric}
	\int_M|\Grad(hv)- \Grad v|_g^2\; \d V = 0. 
\end{equation} 
To do so we observe that for each $i$, the $G$-invariance of $v_i$ guarantees that $\Grad v_i = \Grad(hv_i)$, where equality holds in the sense of $L^2$ vector fields on $M$. 
Since, in addition, $h^{-1}\cdot(\div X)\in L^2(M)\subset [H^1(M)]'$ whenever $X\in \Gamma(TM)$, for any such $X$ we have 
\begin{equation*}
\begin{split}
	\int_Mg(\Grad(h v), X)\; \d V
	& = -\int_Mh v \; \div X\; \d V\\
	& = -\int_Mv \; h^{-1}(\div X)\; \d V\\
	& = -\lim_i\int_Mv_i \; h^{-1}(\div X)\; \d V\\
	& = -\lim_i\int_Mh v_i \; \div X\; \d V\\
	& = \lim_i\int_Mg(\Grad(h v_i), X)\; \d V\\
	& = \lim_i\int_Mg(\Grad v_i, X)\; \d V\\
	& = \int_Mg(\Grad v, X)\; \d V. 
\end{split}
\end{equation*}
 Since $X\in \Gamma(TM)$ is arbitrary, we conclude that $\Grad v = \Grad(h v)$ in the sense of $L^2$ vector fields on $M$.
\end{proof}
For a subgroup $G\subset O(n + 1)$ we denote by $E_G$ the restriction of the functional $E$ given in \eqref{eq:E_functional} to the subspace $H_G^1(\bb S^n)$ of $G$-invariant $H^1(\bb S^n)$ functions. Similarly, if $G$ satisfies property $\mathcal P$ with corresponding element $\tau\in NG\setminus G$ we denote the restriction of $E$ to $\Gamma(G, \tau)$ by $E_\Gamma = E_{\Gamma(G, \tau)}$. In particular, the explicit expression \eqref{eq:E_explicit_expression} holds with $E$ replaced by $E_G$ (respectively $E_{\Gamma(G, \tau)}$) whenever $v\in H_G^1(\bb S^n)$ (respectively $v\in H_{\Gamma(G, \tau)}^1(\bb S^n)$). The following lemma asserts that for suitable $G$, the restriction $E_G$ and the restriction $E_{\Gamma(G, \tau)}$ satisfy the Palais-Smale condition. 
\begin{lemma}
\label{lemma:PS_compactness}
Let $n\geq 3$ and let $G$ be as in \eqref{eq:G_form} for some integer $m\geq 2$ and some integers $n_1, \ldots, n_m$ satisfying both \eqref{eq:nj_sum} and \eqref{eq:nj_largeness}. The restriction $E_G:H_G^1(\bb S^n)\to \bb R$ of the functional $E$ defined in \eqref{eq:E_functional} satisfies the Palais-Smale compactness condition. If in addition to the above hypotheses, $G$ satisfies property $\mathcal P$ with corresponding element $\tau\in NG\setminus G$ then for $\Gamma = \Gamma(G, \tau)$ as in \eqref{eq:Gamma_Gtau}, the restriction $E_\Gamma:H_\Gamma^1(\bb S^n)\to \bb R$ of $E$ satisfies the Palais-Smale compactness condition.
\end{lemma}
\begin{proof}
The proof that $E_\Gamma$ satisfies the Palais-Smale condition whenever $G$ and $\tau$ satisfy property $\mathcal P$ can be obtained by making cosmetic modifications to the proof that $E_G$ satisfies the Palais-Smale condition. For this reason, only a proof that $E_G$ satisfies the Palais-Smale condition will be provided. Accordingly, let $(v_i)_{i = 1}^\infty\subset H_G^1(\bb S^n)$ be a Palais-Smale sequence for $E_G$ so that $E_G(v_i) \to a$ for some $a\in \bb R$ and $E_G'(v_i)\to 0$ in $[H_G^1(\bb S^n)]'$. We emphasize that the convergence $E_G'(v_i)\to 0$ in $[H_G^1(\bb S^n)]'$ means that $\lb E_G'(v_i), \psi\rb \to 0$ for all $\psi\in H_G^1(\bb S^n)$. In particular, on it's own, this convergence does not forbid the existence of $\psi\in H^1(\bb S^n)\setminus H_G^1(\bb S^n)$ for which $\lb E_G'(v_i), \psi\rb \not\to 0$. 
\begin{enumerate}[label = {\bf Step \arabic*.}, wide = 0pt]
	\item We show that $(v_i)_{i = 1}^\infty$ is bounded in $H^1(\bb S^n)$.\\
	For $i$ sufficiently large we have
	\begin{equation*}
	\begin{split}
		\left(1 - \frac{1}{2^*_\mu}\right)\int_{\bb S^n}J_\mu[|v_i|^{2^*_\mu}]|v_i|^{2^*_\mu}\; \d V
		& = 2E_G(v_i) - \lb E_G'(v_i), v_i\rb\\
		& = 2a + \circ(1)(1 + \|v_i\|_{H^1(\bb S^n)})\\
		& \leq 2a + 1 + \circ(1)\|v_i\|_{H^1(\bb S^n)}. 
	\end{split}
	\end{equation*}
	From this estimate and the definition of $E_G$, for $i$ large we have
	\begin{equation*}
	\begin{split}
		\|v_i\|_{H^1(\bb S^n)}^2
		& = 2E_G(v_i) + \frac 1{2^*_\mu}\int_{\bb S^n}J_\mu[|v_i|^{2^*_\mu}]|v_i|^{2^*_\mu}\; \d V\\
		& \leq 2 E_G(v_i) + \frac 1{2^*_\mu - 1}\left(2a + 1 + \circ(1)\|v_i\|_{H^1(\bb S^n)}\right). 
	\end{split}
	\end{equation*}
	Upon rearranging this inequality we obtain (still for $i$ large)
	\begin{equation*}
	\begin{split}
		\|v_i\|_{H^1(\bb S^n)}(\circ(1) + \|v_i\|_{H^1(\bb S^n)})
		& \leq 2E_G(v_i) + \frac1{2^*_\mu - 1}(2a + 1)\\
		& \leq \left(2 + \frac{1}{2^*_\mu - 1}\right)(2|a| + 1), 
	\end{split}
	\end{equation*}
	from which the boundedness of $(v_i)_{i = 1}^\infty$ in $H^1(\bb S^n)$ follows.
	\item We show that there is $v\in H^1(\bb S^n)$ for which $v_i\to v$ in $H^1(\bb S^n)$.\\
	Since $(v_i)_{i = 1}^\infty$ is bounded in $H^1(\bb S^n)$ there is $v\in H^1(\bb S^n)$ and a subsequence of $(v_i)_{i = 1}^\infty$ (whose members are still denoted $v_i$) along which each of the following holds: 
	\begin{equation*}
	\begin{split}
		v_i\rightharpoonup v& \text{ weakly in }H^1(\bb S^n),\\
		v_i\to v & \text{ in }L^2(\bb S^n) \text{ and in }L^{2^*}(\bb S^n)\\
		v_i\to v& \text{ a.e. on }\bb S^n, 
	\end{split}
	\end{equation*}
	where the $L^{2^*}(\bb S^n)$-convergence follows from Corollary \ref{coro:improved_embedding}. For any such $v$, Lemma \ref{lemma:weakly_closed} guarantees that $v\in H_G^1(\bb S^n)$. Moreover, the $L^2(\bb S^n)$ convergence $v_i\to v$ and the weak convergence $v_i\rightharpoonup v$ in $H^1(\bb S^n)$ give
	\begin{equation*}
	\begin{split}
		\|v_i - v\|_{H^1(\bb S^n)}^2
		& = \int_{\bb S^n}\left(|\Grad v_i|^2 -2g(\Grad v_i, \Grad v) + |\Grad v|^2\right)\; \d V + \circ(1)\\
		& = \int_{\bb S^n}\left(|\Grad v_i|^2 - |\Grad v|^2\right)\; \d V + \circ(1).
	\end{split}
	\end{equation*}
	In view of this equality, to show that $v_i\to v$ in $H^1(\bb S^n)$, it suffices to show that 
	\begin{equation}
	\label{eq:sufficient_for_H1_convergence}
		\int_{\bb S^n}|\Grad v_i|^2\; \d V\to \int_{\bb S^n}|\Grad v|^2\; \d V.
	\end{equation} 
	The remainder of the proof is dedicated to doing so. First, since $(v_i)_{i = 1}^\infty$ is bounded in $H^1(\bb S^n)$, since $v\in H_G^1(\bb S^n)$, and since $v_i\to v$ in $L^2(\bb S^n)$ we have 
	\begin{equation}
	\label{eq:EG'vi_vi}
	\begin{split}
		\circ(1)
		& = \circ(1)\|v_i\|_{H^1(\bb S^n)}\\
		& = \lb E_G'(v_i), v_i\rb\\
		& = \int_{\bb S^n}\left(|\Grad v_i|^2 + \frac{n(n - 2)}4 v^2\right)\; \d V - \int_{\bb S^n}J_\mu[|v_i|^{2^*_\mu}]|v_i|^{2^*_\mu}\; \d V + \circ(1). 
	\end{split}
	\end{equation}
	Since in addition, the map $\Phi(\psi)= \int_{\bb S^n} g(\Grad\psi, \Grad v)\; \d V$ satisfies $\Phi\in [H^1(\bb S^n)]'$, using once more the fact that $v_i\rightharpoonup v$ weakly in $H^1(\bb S^n)$ and the fact that $v_i\to v$ in $L^2(\bb S^n)$ we have
	\begin{equation}
	\label{eq:EG'vi_v}
	\begin{split}
		\circ(1)
		& = \circ(1)\|v\|_{H^1(\bb S^n)}\\
		& = \lb E_G'(v_i), v\rb\\
		& = \int_{\bb S^n}\left(g(\Grad v_i, \Grad v)+ \frac{n(n - 2)}4 v_i v\right)\; \d V - \int_{\bb S^n}J_\mu[|v_i|^{2^*_\mu}]|v_i|^{2^*_\mu- 2}v_i v\; \d V\\
		& = \int_{\bb S^n}\left(|\Grad v|^2+ \frac{n(n - 2)}4 v^2\right)\; \d V - \int_{\bb S^n}J_\mu[|v_i|^{2^*_\mu}]|v_i|^{2^*_\mu- 2}v_i v\; \d V+ \circ(1). 
	\end{split}
	\end{equation}
	Combining equations \eqref{eq:EG'vi_vi} and \eqref{eq:EG'vi_v} we have
	\begin{equation}
	\label{eq:now_estimate_error_phi}
	\begin{split}
		\circ(1)
		& = \lb E_G'(v_i), v_i\rb - \lb E_G'(v_i), v\rb\\
		& = \int_{\bb S^n}\left(|\Grad v_i|^2 - |\Grad v|^2\right)\; \d V - \phi_i + \circ(1), 
	\end{split}
	\end{equation}
	where
	\begin{equation*}
		\phi_i = \int_{\bb S^n}J_\mu[|v_i|^{2^*_\mu}]\left(|v_i|^{2^*_\mu} - |v_i|^{2^*_\mu - 2}v_iv\right)\; \d V. 
	\end{equation*}
	Using H\"older's inequality (with exponents $\frac{2n}\mu$, $\frac{2^*}{2^*_\mu - 1}$ and $2^*$), the HLS inequality, and the fact that $v_i\to v$ in $L^{2^*}(\bb S^n)$ we have
	\begin{equation}
	\label{eq:error_phi_bound}
	\begin{split}
		|\phi_i|
		& 
		\leq \int_{\bb S^n}J_\mu[|v_i|^{2^*_\mu}]|v_i|^{2^*_\mu - 1}|v_i - v|\; \d V\\
		& \leq \|J_\mu[|v_i|^{2^*_\mu}]\|_{L^{\frac{2n}\mu}(\bb S^n)}\||v_i|^{2^*_\mu - 1}\|_{L^{\frac{2^*}{2^*_\mu - 1}}(\bb S^n)}\|v_i - v\|_{L^{2^*}(\bb S^n)}\\
		& \leq C\|v_i\|_{L^{2^*}(\bb S^n)}^{2\cdot 2^*_\mu - 1}\|v_i - v\|_{L^{2^*(\bb S^n)}}\\
		& = \circ(1). 
	\end{split}
	\end{equation}
	Bringing \eqref{eq:error_phi_bound} back to \eqref{eq:now_estimate_error_phi} yields \eqref{eq:sufficient_for_H1_convergence}. 
\end{enumerate}

\end{proof}
\section{Existence of Unbounded Sequence of Solutions}
\label{s:existence}
We will use the following version of the Symmetric Mountain Pass Lemma of Ambrosetti and Rabinowitz. See, for example, Theorem 6.5 on p. 114 of \cite{Struwe2008}. 
\begin{oldtheorem}
\label{oldtheorem:SMP}
Let $X$ be an infinite-dimensional Banach space and suppose $\Phi\in C^1(X; \bb R)$ is an even functional that satisfies the Palais-Smale compactness condition and for which $\Phi(0) =0$. Suppose $X = X^-\oplus X^+$, where $X^-$ is finite dimensional and assume the following conditions hold
\begin{enumerate}[label = {\bf \Alph*.}, ref = {\bf \Alph*}, wide = 0pt]
	\item \label{item:1}There exists $\alpha> 0$ and $\rho>0$ such that for all $u\in X^+$ for which $\|u\| = \rho$, the inequality $\Phi(u)\geq\alpha$ holds. 
	\item \label{item:2}For any finite-dimensional subspace $Z\subset X$ there is $R = R(Z)>0$ such that $\Phi(u)\leq 0$ whenever $u\in Z$ and $\|u\|\geq R$. 
\end{enumerate} 
	Then $\Phi$ possesses an unbounded sequence of critical values. 
\end{oldtheorem}
In Proposition \ref{prop:critical_points_restriction} below we apply Theorem \ref{oldtheorem:SMP} to obtain an unbounded sequence of critical points for the restricted functional $E_G$ for suitable $G\subset O(n + 1)$ and to the restricted functional $E_{\Gamma(G, \tau)}$ whenever a suitable $\tau$ exists. In the verification of item \ref{item:2} of Theorem \ref{oldtheorem:SMP} we will use the following lemma, a proof of which can be obtained by making obvious modifications to the proof of Lemma 2.3 in \cite{GaoYang2018}. 
\begin{lemma}
\label{lemma:NL_norm}
The function $\|\cdot\|_{NL(\bb S^n)}: L^{2^*}(\bb S^n)\to \bb R$ defined by 
\begin{equation}
\label{eq:NL_norm}
	\|v\|_{NL(\bb S^n)}^{2\cdot 2^*_\mu}
	= \int_{\bb S^n}J_\mu[|v|^{2^*_\mu}]|v|^{2^*_\mu}\; \d V
\end{equation}
defines a norm on $L^{2^*}(\bb S^n)$. 
\end{lemma}

\begin{prop}
\label{prop:critical_points_restriction}
Let $G$ be of the form \eqref{eq:G_form} for some integer $m\geq 2$ and some integers $n_1, \ldots, n_m$ satisfying both \eqref{eq:nj_sum} and \eqref{eq:nj_largeness}. The restricted functional $E_G:H_G^1(\bb S^n)\to \bb R$ has an unbounded sequence of critical points. If, in addition $G$ and $\tau\in NG\setminus G$ satisfy property $\mathcal P$ then for $\Gamma = \Gamma(G, \tau)$ as in \eqref{eq:Gamma_Gtau}, the restricted functional $E_\Gamma:H_\Gamma^1(\bb S^n)\to \bb R$ has an unbounded sequence of critical points. 
\end{prop}
\begin{proof}
To establish the existence of an unbounded sequence of critical points for $E_G$, we verify that the hypotheses of Theorem \ref{oldtheorem:SMP} hold with $X = H_G^1(\bb S^n)$, with $\Phi = E_G$ and with 
\begin{equation}
\label{eq:X_decompose}
	X^- = X^-(\ell)= \bigoplus_{i = 0}^{\ell - 1}\bb Y_{G, i}
	\qquad \text{ and }\qquad
	X^+ = X^+(\ell) = \overline{\bigoplus_{i = \ell}^\infty\bb Y_{G, i}}
\end{equation}
for suitably large $\ell$, where $\bb Y_{G, i}$ is the space of $G$-invariant spherical harmonics of degree $i$ as defined in \eqref{eq:symmeric_spherical_harmonics}. Evidently $E_G\in C^1(H_G^1(\bb S^n); \bb R)$ is even, satisfies $E_G(0) = 0$, and assumptions \eqref{eq:nj_sum} and \eqref{eq:nj_largeness} ensure that $\dim(H_G^1(\bb S^n)) = +\infty$. Moreover, Lemma \ref{lemma:PS_compactness} guarantees that $E_G$ satisfies the Palais-Smale compactness condition. It remains to establish the existence of $\ell\in \bb N$ for which items \ref{item:1} and \ref{item:2} of Theorem \ref{oldtheorem:SMP} hold with $X^-$ and $X^+$ as in \eqref{eq:X_decompose}. Let $0 = \lambda_0< \lambda_1< \lambda_2< \ldots$ denote the eigenvalues of $-\lap_{\bb S^n}$. First, for any $\ell\in \{1, 2, \ldots\}$ and any $v\in \overline{\bigoplus_{i = \ell}^\infty \bb Y_{G, i}}\setminus\{0\}$ we have
\begin{equation}
\label{eq:ellth_eigenvalue}
\begin{split}
	\frac{\|\Grad v\|_{L^2(\bb S^n)}^2}{\|v\|_{L^2(\bb S^n)}^2}
	& \geq \inf\left\{\frac{\|\Grad w\|_{L^2(\bb S^n)}^2}{\|w\|_{L^2(\bb S^n)}^2}: w\in \overline{\bigoplus_{i = \ell}^\infty \bb Y_{G, i}}\setminus\{0\} \right\}\\
	& \ifdetails
\else \geq
	\fi 
	\inf\left\{\frac{\|\Grad w\|_{L^2(\bb S^n)}^2}{\|w\|_{L^2(\bb S^n)}^2}: w\in \left(\bigoplus_{i = 0}^{\ell - 1} \bb Y_i\right)^\perp\setminus\{0\} \right\}\\
	& = \lambda_\ell, 
\end{split}
\end{equation}
where we used both the containment $\bb Y_{G, i}\subset \bb Y_i$ and the variational characterization of the eigenvalues of $-\lap_{\bb S^n}$.  Corollary \ref{coro:improved_embedding} guarantees the existence of $p>2^*$ for which the estimate in \eqref{eq:improved_embedding} holds. Let us fix any such $p$, set $a = \frac{2(p - 2^*)}{p - 2}\in (0, 2)$ and apply H\"older's inequality with exponents $\frac 2 a$ and $\frac{2}{2 - a}$ to obtain 
\begin{equation*}
	\|v\|_{L^{2^*}(\bb S^n)}^{2^*}
	\leq \|v\|_{L^p(\bb S^n)}^{2^* - a}\|v\|_{L^2(\bb S^n)}^a
	\leq C\|v\|_{H^1(\bb S^n)}^{2^* - a}\|v\|_{L^2(\bb S^n)}^a
\end{equation*}
for all $v\in H_G^1(\bb S^n)$. Using this estimate together with H\"older's inequality and the HLS inequality we find that  
\begin{equation*}
	\int_{\bb S^n}J_\mu[|v|^{2^*_\mu}]|v|^{2^*_\mu}\; \d V
	\leq C\|v\|_{L^{2^*}(\bb S^n)}^{2\cdot 2^*_\mu}
	\leq C\|v\|_{H^1(\bb S^n)}^\sigma\|v\|_{L^2(\bb S^n)}^\delta, 
\end{equation*}
for all $v\in H_G^1(\bb S^n)$, where $\sigma = 2\cdot2^*_\mu(1 - a/2^*)> 0$ and $\delta = 2a\cdot \frac{2^*_\mu}{2^*}> 0$. We note that $\sigma + \delta = 2\cdot 2^*_\mu > 2$. If, in addition, $v\in \overline{\bigoplus_{i = \ell}^\infty\bb Y_{G, i}}\setminus\{0\}$ then continuing in the previous estimate and using estimate \eqref{eq:ellth_eigenvalue} gives
\begin{equation*}
\begin{split}
	\int_{\bb S^n}J_\mu[|v|^{2^*_\mu}]|v|^{2^*_\mu}\; \d V
	& \leq C\|v\|_{H^1(\bb S^n)}^\sigma\|\Grad v\|_{L^2(\bb S^n)}^\delta \left(\frac{\|v\|_{L^2(\bb S^n)}^2}{\|\Grad v\|_{L^2(\bb S^n)}^2}\right)^{\delta/2}\\
	& \leq C\|v\|_{H^1(\bb S^n)}^{\sigma + \delta}\lambda_\ell^{-\delta/2}.
\end{split}
\end{equation*}
Using this estimate and the fact that $\lambda_\ell \to\infty$ as $\ell\to \infty$ we find that if $\ell$ is sufficiently large then for all $v\in \overline{\bigoplus_{i = \ell}^\infty\bb Y_{G, i}}$ for which $\|v\|_{H^1(\bb S^n)} = 1$ we have 
\begin{equation*}
\begin{split}
	E(v)
	&\geq \frac 12\|v\|_{H^1(\bb S^n)}^2 - C\lambda_\ell^{-\delta/2}\|v\|_{H^1(\bb S^n)}^{\sigma + \delta}\\
	& = \left(\frac 12 - C\lambda_\ell^{-\delta/2}\|v\|_{H^1(\bb S^n)}^{\sigma + \delta - 2}\right)\|v\|_{H^1(\bb S^n)}^2\\
	& \geq \frac 14 
\end{split}
\end{equation*}
and thus, for any such $\ell$, item \ref{item:1} of Theorem \ref{oldtheorem:SMP} is satisfied with $\rho = 1$ and $\alpha = 1/4$. To prove that item \ref{item:2} of Theorem \ref{oldtheorem:SMP} holds, let $Z\subset H_G^1(\bb S^n)$ be a finite-dimensional subspace. Since the restriction of the norm $\|\cdot\|_{NL(\bb S^n)}$ (as defined in \eqref{eq:NL_norm}) and the norm $\|\cdot\|_{H^1(\bb S^n)}$ to $Z$ are two norms on a common finite-dimensional subspace there is $C_0= C_0(Z)> 1$ for which 
\begin{equation*}
	C_0^{-1}\|v\|_{NL(\bb S^n)}
	\leq\|v\|_{H^1(\bb S^n)}
	\leq C_0\|v\|_{NL(\bb S^n)}
	\qquad \text{ for all }v\in Z. 
\end{equation*}
Therefore, setting $C_1 = 2^*_\mu C_0^{2\cdot 2^*_\mu}$, if $v\in Z$ satisfies $\|v\|_{H^1(\bb S^n)}\geq (2C_1)^{\frac 1{2(2^*_\mu - 1)}}$ then 
\begin{equation*}
\begin{split}
	E_G(v)
	& \leq \frac 12\|v\|_{H^1(\bb S^n)}^2 - \frac 1{2\cdot2^*_\mu}\left(\frac 1{C_0}\|v\|_{H^1(\bb S^n)}\right)^{2\cdot 2^*_\mu}\\
	& = \frac 12\|v\|_{H^1(\bb S^n)}^2\left(1 - \frac 1{C_1}\|v\|_{H^1(\bb S^n)}^{2\cdot 2^*_\mu - 2}\right)\\
	& \leq -\frac 12\|v\|_{H^1(\bb S^n)}^2\\
	& \leq 0, 
\end{split}
\end{equation*}
so item \ref{item:2} of Theorem \ref{oldtheorem:SMP} is satisfied. Applying Theorem \ref{oldtheorem:SMP} guarantees the existence of a sequence of critical points $(v_i)_{i = 1}^\infty\subset H_G^1(\bb S^n)$ of $E_G$ for which $E_G(v_i)\to \infty$. Since the inequality $2E_G(v)\leq \|v\|_{H^1(\bb S^n)}$ holds for all $v\in H^1(\bb S^n)$, any such sequence $(v_i)_{i = 1}^\infty$, satisfies $\|v_i\|_{H^1(\bb S^n)}\to \infty$. \\
In the case that $G$ and $\tau$ satisfy property $\mathcal P$ the proof that there is an unbounded sequence of critical points for $E_\Gamma$ is essentially the same as the proof that $E_G$ has an unbounded sequence of critical points. The only difference is that in this case the existence of $\xi\in \bb S^n$ for which $\tau\xi\not\in G\xi$ is used to ensure that $\dim(H_\Gamma^1(\bb S^n)) = +\infty$. 
\end{proof}
\begin{proof}[Proof of Theorem \ref{theorem:main}]
To facilitate the discussion in Subsection \ref{ss:remarks_and_examples} below, the proof is divided into two separate but similar constructions.
\begin{enumerate}[label = {\bf Construction \arabic*.}, ref = {\bf Construction \arabic*}, wide = 0pt]
	\item \label{item:construction1} Let $G\subset O(n + 1)$ be of the form \eqref{eq:G_form} for some integer $m\geq 2$ and some integers $n_1, \ldots, n_m$ satisfying both \eqref{eq:nj_sum} and \eqref{eq:nj_largeness}. Proposition \ref{prop:critical_points_restriction} ensures the existence of a sequence $(v_i)_{i= 1}^\infty\subset H_G^1(\bb S^n)$ of critical points for $E_G$ for which $\|v_i\|_{H^1(\bb S^n)}\to \infty$ as $i\to\infty$. Moreover, according to Lemma \ref{lemma:E_invariant}, the energy functional $E$ is invariant under the action of $G$ so the Principle of Symmetric Criticality \cite{Palais1979} guarantees that each $v_i$ is a critical point of $E$. With $P$ as in \eqref{eq:P_projection}, Proposition \ref{prop:critical_point_correspondence} guarantees that $u_i = Pv_i$ is a critical point of $F$ and Lemma \ref{lemma:extend_P} guarantees that $\|u_i\|_{D^{1, 2}(\bb R^n)} = \|v_i\|_{H^1(\bb S^n)}\to \infty$. If either $n\in \{3, 4\}$ or if both $n\geq 5$ and $\mu\in (0, 4]$ then all but finitely many of the $u_i$ are sign-changing. Indeed, suppose $i$ is an index for which $u_i$ is not sign-changing and thus either $u_i$ or $-u_i$ is a positive solution to problem \eqref{eq:choquard}. The classification theorem of \cite{DuYang2019} implies that either $u_i$ or $-u_i$ is of the form \eqref{eq:scaled_translated_bubble} for some $(x_0, \lambda)\in \bb R^n\times (0, \infty)$. Since $\|U_{x_0, \lambda}\|_{D^{1, 2}(\bb R^n)}$ is independent of $(x_0, \lambda)$, for any such $i$ we have $\|u_i\|_{D^{1, 2}(\bb R^n)} = \|U\|_{D^{1, 2}(\bb R^n)}$ where $U$ is as in \eqref{eq:UnscaledBubble}. Therefore, if $i_0\in \bb N$ is large enough so that $\|u_j\|_{D^{1, 2}(\bb R^n)}> 1 + \|U\|_{D^{1, 2}(\bb R^n)}$ whenever $j\geq i_0$ then $i$ satisfies $i\leq i_0$. 
\item \label{item:construction2} In the case that $G$ can be chosen to satisfy property $\mathcal P$ with corresponding element $\tau\in NG\setminus G$ (in view of Example \ref{example:standard}, this holds whenever $n\neq 4$), arguing as above we obtain an unbounded sequence of critical points $(v_i)_{i = 1}^\infty$ of $E$ for which $(v_i)_{i =1}^\infty\subset H_{\Gamma(G, \tau)}^1(\bb S^n)$. For each $i$ for which $\|v_i\|_{H^1(\bb S^n)}\neq 0$ (in particular, for infinitely many $i$) the invariance of $v_i$ under the action $\beta$ defined in \eqref{eq:non_standard_action} implies that $v_i\circ\tau^{-1} = -v_i$, so $v_i$ changes sign. With $P$ as in \eqref{eq:P_projection} we find that $Pv_i$ is a sign-changing solution to \eqref{eq:choquard}. 
\end{enumerate}
\end{proof}
\subsection{Remarks and Examples}
\label{ss:remarks_and_examples}
In this subsection we record some remarks and give some explicit examples of the two constructions in the proof of Theorem \ref{theorem:main}. 
\begin{remark}
\label{remark:nodal_info}
\ref{item:construction2} in the proof of Theorem \ref{theorem:main} yields information on the nodal sets of the solutions. To see this, suppose $G\subset O(n + 1)$ satisfies property $\mathcal P$ with corresponding element $\tau\in NG\setminus G$. If $v\in H^1_{\Gamma(G, \tau)}(\bb S^n)$ (in particular, for each $v_i$ in \ref{item:construction2} in the proof of Theorem \ref{theorem:main}), then $\{\xi\in \bb S^n: \tau\xi \in G\xi\}\subset\{\xi\in \bb S^n: v(\xi) = 0\}$. Indeed, if $\xi\in \bb S^n$ and if $h\in G$ with $\tau\xi = h\xi$ then
\begin{equation*}
	v(\xi) = v(h\xi) = v(\tau\xi) = -v(\xi). 
\end{equation*}
\end{remark}
\begin{remark}
A direct computation using the explicit expression of $E$ in \eqref{eq:E_explicit_expression} shows that
\begin{equation*}
	\lb E'(h\cdot v), \varphi\rb
	= \lb E'(v), h^{-1}\cdot \varphi\rb
\end{equation*}	
whenever $v, \varphi\in H^1(\bb S^n)$ and $h\in O(n + 1)$. In particular, each critical point of $E$ generates an $O(n+ 1)$-orbit of critical points of $E$. 
\end{remark}
In the following examples we discus explicit cases of both constructions in the proof of Theorem \ref{theorem:main} (when applicable) in the cases $n = 3, 4, 5$. 
\begin{example}
If $n = 3$ then the only subgroup $G\subset O(4)$ of the form \eqref{eq:G_form} with $m\geq 2$ that satisfies both \eqref{eq:nj_sum} and \eqref{eq:nj_largeness} is $G = O(2)\times O(2)$. \ref{item:construction1} yields critical points $v_i$ of $E$ that depend on $\xi = (\xi_1, \xi_2)\in \bb S^3\subset \bb R^2\times \bb R^2$ through only $|\xi_1|$ and $|\xi_2|$. All but finitely many of these critical points are sign-changing independently of the value of $\mu\in (0, 3)$. Since $G$ satisfies property $\mathcal P$ with 
\begin{equation*}
	\tau = \begin{bmatrix}
	0 & I_2\\
	I_2 & 0
	\end{bmatrix}
	\in O(4)\setminus G, 
\end{equation*}
\ref{item:construction2} yields critical points $v_i$ of $E$ that may be the same as those obtained from \ref{item:construction1}. However \ref{item:construction2} yields the nodal set information $\{\xi\in \bb S^3: \tau\xi\in G\xi\}\subset\{\xi\in \bb S^3: v_i(\xi)= 0\}$ as in Remark \ref{remark:nodal_info}. 
\end{example}
\begin{example}
In the case $n = 4$ we do not have an example of a subgroup $G\subset O(5)$ together with an element $\tau\in NG\setminus G$ for which property $\mathcal P$ holds, so we do not know whether \ref{item:construction2} can be applied. On the other hand, the group $G= O(3)\times O(2)$ is the only subgroup of $O(5)$ (up to permutation of coordinates) of the form \eqref{eq:G_form} for some $m\geq 2$ and some integers $n_1, \ldots, n_m$ that satisfy both \eqref{eq:nj_sum} and \eqref{eq:nj_largeness}. For this $G$, \ref{item:construction1} yields an $H^1(\bb S^n)$-unbounded sequence of critical points $(v_i)_{i = 1}^\infty$ for $E$ that depend on $\xi = (\xi_1, \xi_2)\in \bb S^4\subset \bb R^3\times \bb R^2$ through only $|\xi_1|$ and $|\xi_2|$. All but finitely many of these functions are sign-changing independently of the value of $\mu\in (0, 4)$.  
\end{example}
\begin{example}
In the case $n = 5$ there are three subgroups of $O(6)$ of the form \eqref{eq:G_form} for some $m\geq 2$ and some integers $n_1, \ldots, n_m$ that satisfy both \eqref{eq:nj_sum} and \eqref{eq:nj_largeness}. These subgroups are (up to permutation of coordinates)
\begin{equation*}
\begin{split}
	G_1 & = O(4)\times O(2)\\
	G_2 & = O(2)\times O(2)\times O(2)\\
	G_3 & = O(3)\times O(3).
\end{split}
\end{equation*}
In what follows we use the notation $\xi = (\xi_1, \ldots, \xi_6)\in \bb S^5$. \\
Applying \ref{item:construction1} with $G = G_1$ yields an $H^1(\bb S^n)$-unbounded sequence of critical points $(v_i^{(1)})_{i = 1}^\infty\subset H_{G_1}^1(\bb S^n)$ for $E$ that satisfy
\begin{equation*}
	v_i^{(1)}(\xi) = v_i^{(1)}(\xi_1^2 + \ldots +\xi_4^2, \xi_5^2 + \xi_6^2).  
\end{equation*}
If $\mu\in (0, 4]$ then all but finitely many of these functions are sign-changing. We do not know whether this property also holds for $\mu\in (4, 5)$. \\
Applying \ref{item:construction1} with $G = G_2$ yields an $H^1(\bb S^n)$-unbounded sequence of critical points $(v_i^{(2)})_{i = 1}^\infty\subset H_{G_2}^1(\bb S^n)$ for $E$ that satisfy
\begin{equation*}
	v_i^{(2)}(\xi) = v_i^{(2)}(\xi_1^2 + \xi_2^2, \xi_3^2 + \xi_4^2, \xi_5^2 + \xi_6^2).  
\end{equation*}
If $\mu\in (0, 4]$ then all but finitely many of these functions are sign-changing. We do not know whether this property also holds for $\mu\in (4, 5)$. Observe that $G_2\subset G_1$ and thus $H_{G_1}^1(\bb S^n)\subset H_{G_2}^1(\bb S^n)$. In particular, it is possible that the containment $\{v_i^{(2)}\}_{i = 1}^\infty \subset \{v_i^{(1)}\}_{i = 1}^\infty$ holds. For 
\begin{equation*}
	\tau_2 = \begin{bmatrix}
		0 & I_2 & 0\\
		I_2 & 0 & 0\\
		0 & 0 & I_2
	\end{bmatrix}
	\in O(6)\setminus G_2, 
\end{equation*}
applying \ref{item:construction2} with $\Gamma(G_2, \tau_2)$ yields an $H^1(\bb S^n)$-unbounded sequence $(w_i^{(2)})_{i = 1}^\infty\subset H_{\Gamma(G_2, \tau_2)}^1(\bb S^n)\setminus H_{G_1}^1(\bb S^n)\subset H_{G_2}^1(\bb S^n)\setminus H_{G_1}^1(\bb S^n)$ of critical points of $E$. These functions satisfy
\begin{equation*}
	w_i^{(2)}(\xi_1^2 + \xi_2^2, \xi_3^2 + \xi_4^2, \xi_5^2 + \xi_6^2)
	= w_i^{(2)}(\xi)
	= -w_i^{(2)}(\xi_3^2 + \xi_4^2, \xi_1^2 + \xi_2^2, \xi_5^2 + \xi_6^2)
\end{equation*}
for all $\xi\in \bb S^n$ and are sign-changing independently of the value of $\mu\in (0, 5)$. \\
Applying \ref{item:construction1} with $G = G_3$ yields an $H^1(\bb S^n)$-unbounded sequence of critical points $(v_i^{(3)})_{i = 1}^\infty\subset H_{G_3}^1(\bb S^n)$ for $E$ that satisfy 
\begin{equation*}
	v_i^{(3)}(\xi) = v_i^{(3)}(\xi_1^2 + \xi_2^2 + \xi_3^2, \xi_4^2 + \xi_5^2 + \xi_6^3). 
\end{equation*}
If $\mu\in (0, 4]$ then all but finitely many of these functions are sign-changing. We do not know whether this property also holds when $\mu\in (4, 5)$. For 
\begin{equation*}
	\tau_3 = \begin{bmatrix}
	0 & I_3\\
	I_3 & 0
	\end{bmatrix}
	\in O(6)\setminus G_3
\end{equation*}
applying \ref{item:construction2} with $\Gamma(G_3, \tau_3)$ yields an $H^1(\bb S^n)$-unbounded sequence $(w_i^{(3)})_{i = 1}^\infty\subset H_{\Gamma(G_3, \tau_3)}^1(\bb S^n)$ of critical points of $E$. These functions satisfy
\begin{equation*}
	w_i^{(3)}(\xi_1^2 + \xi_2^2 + \xi_3^2, \xi_4^2 + \xi_5^2 + \xi_6^2)
	= w_i^{(3)}(\xi)
	- w_i^{(3)}(\xi_4^2 + \xi_5^2 + \xi_6^2, \xi_1^2 + \xi_2^2 + \xi_3^2)
\end{equation*}
and are sign-changing independently of the value of $\mu\in (0, 5)$. 
\end{example}
%
\appendix\section{Regularity of solutions}
\label{s:appendix}
%
%
In this section we provide an outline of a proof of Theorem \ref{theorem:regularity}. For ease of notation, throughout this section we use the notational convention $\|u\|_p = \|u\|_{L^p(\bb R^n)}$. The assertion of Theorem \ref{theorem:regularity} follows from the following proposition. 
\begin{prop}
\label{prop:system_regularity}
Let $n\geq 3$, let $\mu\in (0, n)$ and let $2^*_\mu$ be as in \eqref{eq:the_exponent}. There exists $\alpha\in (0, 1)$ such that if $(u, v)\in L^{2^*}(\bb R^n)\times L^{\frac{2n}\mu}(\bb R^n)$ satisfies
\begin{equation}
\label{eq:the_system}
\left\{\begin{split}
	u(x) & = \int_{\bb R^n}\frac{v(y)|u(y)|^{2^*_\mu - 2}u(y)}{|x - y|^{n - 2}}\; \d y\\
	v(x) & = I_\mu[|u|^{2^*_\mu}](x),
\end{split}
\right.
\end{equation}
then $u\in C^{2, \alpha}(\bb R^n)$. 
\end{prop}
The primary technical challenge in the proof of Proposition \ref{prop:system_regularity} is to establish the local boundedness of solutions to \eqref{eq:the_system}. This challenge is addressed in the following lemma. 
\begin{lemma}
\label{lemma:local_boundedness}
Let $n\geq 3$ and let $\mu\in (0, n)$. If $(u,v)\in L^{2^*}(\bb R^n)\times L^{\frac{2n}\mu}(\bb R^n)$ is a solution to \eqref{eq:the_system} then $u\in L^\infty_{\loc}(\bb R^n)$ and $v\in L^\infty_{\loc}(\bb R^n)$. 
\end{lemma}
 The following integrability boosting lemma will be used in the proof of Lemma \ref{lemma:local_boundedness}. 
\begin{lemma}
\label{lemma:dual_integrability_boosting}
Let $n\geq 3$, let $\mu\in (0, n)$, and suppose $q$ satisfies
\begin{equation}
\label{eq:dual_q_range}
	\frac{1}{2^*}- \frac\mu{2n} < \frac 1q < \frac1{2^*}. 
\end{equation}
For $R>0$, if $h, w\in L^{2^*}(B_R)$ and $U\in L^{\frac{2n}{n + 2 - \mu}}(B_R)$ are non-negative functions for which both $h|_{B_{R/2}}\in L^q(B_{R/2})$ and 
\begin{equation}
\label{eq:boost_main_inequality_assumption}
	w(x)
	\leq \int_{B_R}\frac{U(y)}{|x - y|^{n-2}}\int_{B_R}\frac{U(z)w(z)}{|y - z|^\mu}\; \d z \; \d y + h(x)
	\qquad \text{ for all }x\in B_R
\end{equation}
then there is $\epsilon = \epsilon(n, \mu, q)>0$ and $C = C(n, \mu, q,  \epsilon)>0$ such that if $\|U\|_{L^{\frac{2n}{n+ 2 - \mu}}(B_R)}< \epsilon$ then 
\begin{equation}
\label{eq:boost_asserted_estimate}
	\|w\|_{L^q(B_{R/4})}
	\leq C\left(R^{n\left(\frac 1q - \frac{1}{2^*}\right)}\|w\|_{L^{2^*}(B_R)} + \|h\|_{L^q(B_{R/2})}\right). 
\end{equation}
\end{lemma}
Up to cosmetic modifications, the proof of Lemma \ref{lemma:dual_integrability_boosting} is the same as the proof of Proposition 4.4 of \cite{DouZhu2015} so we omit the details. \\

For any solution $(u,v)\in L^{2^*}(\bb R^n)\times L^{\frac{2n}\mu}(\bb R^n)$ to \eqref{eq:the_system} and for any $0< r< R$ we have the pointwise estimate
\begin{equation*}
	|u(x)|
	\leq \int_{B_R}\frac{|u(y)|^{2^*_\mu - 1}}{|x- y|^{n - 2}}\int_{B_R}\frac{|u(z)|^{2^*_\mu - 1}|u(z)|}{|y - z|^\mu}\; \d z \; \d y + h_{r, R}(x),
\end{equation*}
where 
\begin{equation}
\label{eq:h_rR}
\begin{split}
	h_{r, R}(x)
	= & \; \int_{B_r}\frac{|u(y)|^{2^*_\mu - 1}}{|x- y|^{n - 2}}\int_{\bb R^n\setminus B_R}\frac{|u(z)|^{2^*_\mu - 1}|u(z)|}{|y - z|^\mu}\; \d z\; \d y\\
	& + \int_{\bb R^n\setminus B_r}\frac{v(y)|u(y)|^{2^*_\mu - 1}}{|x - y|^{n - 2}}\; \d y. 
\end{split}
\end{equation}
The following lemma verifies the integrability assumption for $h_{r, R}$ necessary for applying Lemma \ref{lemma:dual_integrability_boosting}. 
\begin{lemma}
\label{lemma:h_integrability}
Let $n\geq 3$, let $\mu\in (0, n)$, let $(u,v)\in L^{2^*}(\bb R^n)\times L^{\frac{2n}\mu}(\bb R^n)$ be a solution to \eqref{eq:the_system} and let $h_{r, R}$ be as in \eqref{eq:h_rR}. If  $0< r<R< \frac{3r}2$ then $h_{r, R}\in L^{2^*}(B_R)$ and the following containments and estimates hold for $h_{r, R}$: 
\begin{enumerate}[label = {\bf(\alph*)}, wide = 0pt]
	\item If $\mu\in (0, n -2)$ then $h_{r, R}\in L^{\frac{2n}{n - 2 - \mu}}(B_{R/2})$ and there is a constant $C(n, \mu)>0$ such that 
\begin{equation}
\label{eq:hrR_Lq_estimate_musmall}
	\|h_{r, R}\|_{L^{\frac{2n}{n - 2 - \mu}}(B_{R/2})}
	\leq \frac{C\|u\|_{2^*}^{2^*_\mu - 1}}{R^{\frac \mu 2}}\left(\|v\|_{\frac{2n}\mu} + \frac{R^\mu}{(R - r)^\mu}\|u\|_{2^*}^{2^*_\mu}\right).
\end{equation}
	\item If $\mu = n - 2$ then for any $q\in (2^*, \infty)$ we have $h_{r, R}\in L^q(B_{R/2})$ and there is a constant $C(n, q)>0$ such that 
	\begin{equation}
	\label{eq:hrR_Lq_estimate_mumedium}
		\|h_{r,R}\|_{L^q(B_{R/2})}
		\leq \frac{C\|u\|_{2^*}^{\frac 4{n - 2}}}{R^{n\left(\frac{1}{2^*}- \frac 1q \right)}}\left(\|v\|_{2^*} + \frac{R^{n - 2}}{(R - r)^{n - 2}}\|u\|_{2^*}^{\frac{n + 2}{n - 2}}\right). 
	\end{equation}
	\item If $\mu\in (n - 2, n)$ then $h_{r, R}\in L^\infty(B_{R/2})$ and there is a constant $C(n, \mu)>0$ such that
	\begin{equation}
	\label{eq:hrR_Linfty_estimate_mularge}
		\|h_{r, R}\|_{L^\infty(B_{R/2})}
		\leq \frac{C\|u\|_{2^*}^{2^*_\mu - 1}}{R^{\frac{n - 2}2}}\left(\|v\|_{\frac{2n}{\mu}} + \frac{R^\mu}{(R - r)^\mu}\|u\|_{2^*}^{2^*_\mu}\right). 
	\end{equation}
\end{enumerate}
\end{lemma}
\begin{proof}
For $0< r< R$ define 
\begin{equation}
\label{eq:ur}
	u_r(x)
	= \int_{\bb R^n\setminus B_r}\frac{v(y)|u(y)|^{2^*_\mu - 1}}{|x- y|^{n - 2}}\; \d y, 
\end{equation}
so that
\begin{equation}
\label{eq:hrR_decomposition}
	h_{r, R}(x)
	= \int_{B_r}\frac{|u(y)|^{2^*_\mu - 1}}{|x- y|^{n - 2}}\int_{\bb R^n\setminus B_R}\frac{|u(z)|^{2^*_\mu}}{|y- z|^\mu}\; \d z\; \d y + u_r(x). 
\end{equation}
We assume henceforth that $0< r< R< \frac{3r}{2}$. For any $x\in B_{R/2}$ and $y\in \bb R^n\setminus B_r$ we have
$
	|x- y|
	\geq \frac{|y|}4,  
$
so for any such $x$, applying H\"older's inequality 
gives
\begin{equation}
\label{eq:ur_Linfty_estimate}
\begin{split}
	|u_r(x)|
	& \leq C(n)\int_{\bb R^n\setminus B_r}\frac{v(y)|u(y)|^{2^*_\mu - 1}}{|y|^{n - 2}}\; \d y\\
	& \leq C(n)\|v\|_{\frac{2n}\mu}\||u|^{2^*_\mu - 1}\|_{\frac{2^*}{2^*_\mu - 1}}\left(\int_{\bb R^n\setminus B_r}|y|^{-2n}\right)^{\frac{1}{2^*}}\\
	& \leq C(n)r^{-\frac{n - 2}2}\|v\|_{\frac{2n}{\mu}}\|u\|_{2^*}^{2^*_\mu - 1}. 
\end{split}
\end{equation}
Using this pointwise estimate, we find that for any $q\in [1, \infty)$ there is $C = C(n, q)>0$ such that
\begin{equation}
\label{eq:ur_Lq_estimate}
\begin{split}
	\|u_r\|_{L^q(B_{R/2})}
	& \leq Cr^{-\frac{n - 2}2}R^{\frac nq}\|v\|_{\frac{2n}{\mu}}\|u\|_{2^*}^{2^*_\mu - 1}\\
	& \leq CR^{n\left(\frac 1q- \frac1{2^*}\right)}\|v\|_{\frac{2n}{\mu}}\|u\|_{2^*}^{2^*_\mu - 1}.
\end{split}
\end{equation}
Independently, we have an $L^\infty(B_r)$ estimate for $I_\mu[|u|^{2^*_\mu}\chi_{\bb R^n\setminus B_R}]$. Indeed, setting $\theta = \frac{r}{R}\in (0, 1)$, for any $z\in \bb R^n\setminus B_R$ and $y\in B_r$ we have
$
	|y - z|
	\geq(1 - \theta)|z|, 
$
so (still for $y\in B_r$) H\"older's inequality gives
\begin{equation}
\label{eq:Imu_tail_Linfty_estimate}
\begin{split}
	\int_{\bb R^n\setminus B_R}\frac{|u(z)|^{2^*_\mu}}{|y - z|^\mu}\; \d z
	& \leq \frac{1}{(1-\theta)^\mu}\int_{\bb R^n\setminus B_R}\frac{|u(z)|^{2^*_\mu}}{|z|^\mu}\;\d z\\
	& \leq \frac{1}{(1 - \theta)^\mu}\|u\|_{2^*}^{2^*_\mu}\left(\int_{\bb R^n\setminus B_R}|z|^{-2n}\; \d z\right)^{1 - \frac{2^*_\mu}{2^*}}\\
	& \leq \frac{C(n,\mu)R^{\frac\mu 2}}{(R - r)^\mu }\|u\|_{2^*}^{2^*_\mu}. 
\end{split}
\end{equation}
In view of this estimate and \eqref{eq:hrR_decomposition} we have
\begin{equation}
\label{eq:hrR_ur_pointwise_estimate}
	|h_{r, R}(x) - u_r(x)|
	\leq \frac{C(n,\mu)R^{\frac\mu 2}}{(R - r)^\mu }\|u\|_{2^*}^{2^*_\mu}\int_{B_r}\frac{|u(y)|^{2^*_\mu - 1}}{|x- y|^{n - 2}}\; \d y
\end{equation}
for $x\in B_{R/2}$. Since the estimate of 
$
	x\mapsto \int_{B_r}\frac{|u(y)|^{2^*_\mu - 1}}{|x- y|^{n - 2}}\;\d y
$
for $x\in B_{R/2}$ depends on whether $\mu\in (0, n- 2)$, $\mu = n - 2$ or $\mu\in (n - 2, n)$, the remainder of the proof considers these cases separately. 
\begin{enumerate}[label = {\bf Case \arabic*.}, ref = {\bf Case \arabic*}, wide = 0pt]
	\item 
	Suppose $\mu\in (0, n - 2)$. For $q = \frac{2n}{n - 2 - \mu}\in (2^*, \infty)$ we have $\frac1q = \frac{2^*_\mu - 1}{2^*} - \frac 2n$ so the HLS inequality gives
	\begin{equation*}
	\begin{split}
		\left\|\int_{B_r}\frac{|u(y)|^{2^*_\mu - 1}}{|\cdot - y|^{n - 2}}\; \d y\right\|_{L^q(B_{R/2})}
		& \leq C(n, \mu)\||u|^{2^*_\mu - 1}\|_{L^{\frac{2^*}{2^*_\mu - 1}}(B_r)}\\
		& \leq C(n, \mu)\|u\|_{L^{2^*}(\bb R^n)}^{2^*_\mu - 1}. 
	\end{split}
	\end{equation*}
	Combining this with \eqref{eq:hrR_ur_pointwise_estimate} yields
	\begin{equation*}
		\|h_{r, R}- u_r\|_{L^{\frac{2n}{n - 2 - \mu}}(B_{R/2})}
		\leq \frac{C(n, \mu)R^{\frac\mu 2}}{(R - r)^\mu}\|u\|_{L^{2^*}(\bb R^n)}^{2\cdot 2^*_\mu - 1}. 
	\end{equation*}
	Combining this estimate and estimate \eqref{eq:ur_Lq_estimate} (applied with $q = \frac{2n}{n - 2 - \mu}$) gives estimate \eqref{eq:hrR_Lq_estimate_musmall}. 
	\item 
	Suppose $\mu = n - 2$. For any $q\in (2^*, \infty)$ define $p$ by $\frac 1p = \frac 1q + \frac 2n\in \left(\frac 2n, \frac{n + 2}{2n}\right) $ so that $|u|^{2^*_\mu - 1} = |u|^{\frac 4{n - 2}}\in L^{\frac n2}(B_r)\subset L^p(B_r)$ and
	\begin{equation*}
	\begin{split}
		\||u|^{2^*_\mu - 1}\|_{L^p(B_r)}^p
		& = \int_{B_r}|u|^{\frac{4p}{n - 2}}\\
		& \leq Cr^{n - 2p}\|u\|_{2^*}^{\frac{4p}{n - 2}} 
	\end{split}
	\end{equation*}
	for some constant $C = C(n, q)>0$. Combining this estimate with the HLS inequality gives
	\begin{equation*}
	\begin{split}
		\left\|\int_{B_r}\frac{|u(y)|^{2^*_\mu - 1}}{|\cdot - y|^{n- 2}}\; \d y\right\|_{L^q(B_{R/2})}
		& \leq C\||u|^{2^*_\mu - 1}\|_{L^p(B_r)}\\
		& \ifdetails
\else \leq
		\fi 
		Cr^{\frac nq}\|u\|_{2^*}^{\frac4{n- 2}},  
	\end{split}
	\end{equation*}
	where $C = C(n, q)$. Combining this with estimate \eqref{eq:hrR_ur_pointwise_estimate} (applied with $\mu = n -2$) we obtain 
	\begin{equation*}
	\begin{split}
		\|h_{r, R}- u_r\|_{L^q(B_{R/2})}
		& \leq 
		\frac{C(n,q)R^{\frac{n -2} 2}r^{\frac nq}}{(R - r)^{n - 2}}\|u\|_{2^*}^{\frac{n + 6}{n - 2}}\\
		& \leq \frac{C(n,q)R^{n\left(\frac 1q + \frac 1{2*}\right)}}{(R - r)^{n - 2}}\|u\|_{2^*}^{\frac{n + 6}{n - 2}}.
	\end{split}
	\end{equation*}
	Combining this estimate with estimate \eqref{eq:ur_Lq_estimate} (applied with $\mu = n - 2$) gives estimate \eqref{eq:hrR_Lq_estimate_mumedium}. 
	\item 
	Suppose $\mu\in (n - 2, n)$. In this case we have $\frac{2n(n - 2)}{n - 2 + \mu}< n$ so for any $x\in B_{R/2}$, H\"older's inequality gives
	\begin{equation*}
	\begin{split}
		\int_{B_r}\frac{|u(y)|^{2^*_\mu- 1}}{|x- y|^{n - 2}}\;\d y
		& 
		\leq
		\|u\|_{2^*}^{2^*_\mu - 1}\left(\int_{B_{2R}(x)}|x- y|^{-\frac{2n(n -2)}{n - 2 + \mu}}\; \d y\right)^{1 - \frac{2^*_\mu- 1}{2^*}}\\
		& \leq C(n, \mu)R^{\frac{2 + \mu - n}{2}}\|u\|_{2^*}^{2^*_\mu- 1}. 
	\end{split}
	\end{equation*}
	Bringing this estimate back to \eqref{eq:hrR_ur_pointwise_estimate} gives
	\begin{equation*}
		\|h_{r, R}- u_r\|_{L^\infty(B_{R/2})}
		\leq \frac{C(n, \mu)R^{1 + \mu - \frac n2}}{(R - r)^\mu}\|u\|_{L^{2^*}(\bb R^n)}^{2\cdot 2^*_\mu - 1},
	\end{equation*}
	and combining this estimate with \eqref{eq:ur_Linfty_estimate} gives \eqref{eq:hrR_Linfty_estimate_mularge}. 
\end{enumerate}
\end{proof}
The following lemma shows that every solution $(u,v)$ to \eqref{eq:the_system} is bounded in some neighborhood of the origin. 
\begin{lemma}
\label{lemma:bounded_near_origin}
Let $n\geq 3$ and let $\mu\in (0, n)$. If $(u,v)\in L^{2^*}(\bb R^n)\times L^{\frac{2n}\mu}(\bb R^n)$ is a solution to \eqref{eq:the_system} then there is $\sigma>0$ for which $u, v\in L^\infty(B_\sigma)$. 
\end{lemma}
\begin{proof}
We separately consider the case $\mu\in [n - 2, n)$ and the case $\mu\in (0, n - 2)$. 
\begin{enumerate}[label = {\bf Case \arabic*.}, ref = {\bf Case \arabic*}, wide = 0pt]
	\item \label{item:origin_bound_case1} Assume $\mu\in [n - 2, n)$. 
	The inequalities
	\begin{equation*}
		\frac 1{2^*}-\frac\mu{2n}< \frac{n- \mu}n\cdot\frac{1}{2^*_\mu}< \frac 1{2^*}
	\end{equation*}
	hold independently of whether $\mu\in [n - 2, n)$ or $\mu\in (0, n-2)$ 
	so we fix $q$ for which 
	\begin{equation}
	\label{eq:case1_q}
		\frac{1}{2^*} - \frac\mu{2n}
		< \frac 1q
		< \frac{n - \mu}n\cdot\frac{1}{2^*_\mu}. 
	\end{equation}
	For such $q$ we have $q>2^*$ and the assumption $\mu\geq n - 2$ guarantees that
	\begin{equation}
	\label{eq:case1_q_large}
		\frac 1q< \frac 2n\cdot\frac{1}{2^*_\mu - 1}. 
	\end{equation}
	For any positive numbers $r$ and $R$ satisfying $0< r< R< \frac{3r}2$, Lemma \ref{lemma:h_integrability} guarantees that the function $h_{r, R}$ defined in \eqref{eq:hrR_decomposition} satisfies $h_{r, R}\in L^{2^*}(B_R)\cap L^q(B_{R/2})$. Applying Lemma \ref{lemma:dual_integrability_boosting} with $w = |u|$ and $U = |u|^{2^*_\mu -1}$ guarantees the existence of $\epsilon = \epsilon(n, \mu, q)>0$ (independent of $R$) for which $u\in L^q(B_{R/4})$ whenever $\|u\|_{L^{2^*}(B_R)}^{2^*_\mu- 1}< \epsilon$. Fix any such $\epsilon$ then choose $R = R(\epsilon)>0$ sufficiently small so that $\|u\|_{L^{2^*}(B_{4R})}^{2^*_\mu - 1}< \epsilon$ and hence $u\in L^q(B_R)$. Now decompose $v(x)$ as
	\begin{equation}
	\label{eq:v_Br4_decomposition}
		v(x)
		= \int_{B_R}\frac{|u(y)|^{2^*_\mu}}{|x-  y|^\mu}\; \d y + v_R(x), 
	\end{equation}	
	where, for $\rho>0$ we set
	\begin{equation*}
		v_{\rho}(x):= \int_{\bb R^n\setminus B_\rho}\frac{|u(y)|^{2^*_\mu}}{|x-  y|^\mu}\; \d y. 
	\end{equation*}
	A routine argument using H\"older's inequality and the assumption $u\in L^{2^*}(\bb R^n)$ shows that $v_R\in L^\infty_{\loc}(B_R)$. Moreover, for any $x\in B_R$, H\"older's inequality and assumption \eqref{eq:case1_q} give
	\begin{equation}
	\label{eq:case1_v-vR_pointwise}
	\begin{split}
		|v(x) - v_R(x)|
		& \leq \|u\|_{L^q(B_R)}^{2^*_\mu}\left(\int_{B_{2R}(x)}|x- y|^{-\frac{q\mu}{q - 2^*_\mu}}\; \d y\right)^{1 - \frac{2^*_\mu}q}\\
		& \leq C(n, \mu, q)R^{n(\frac{n- \mu}n - \frac{2^*_\mu}q)}\|u\|_{L^q(B_R)}^{2^*_\mu}.
	\end{split}
	\end{equation}
	Combining this estimate with the fact that $v_R\in L^\infty_{\loc}(B_R)$ shows that $v\in L^\infty_{\loc}(B_R)$. Now using the fact that $v\in L^\infty_{\loc}(B_R)$ we obtain the pointwise estiamte
	\begin{equation}
	\label{eq:case1_absu_pointwise}
	\begin{split}
		|u(x)|
		& \leq\int_{B_{R/2}}\frac{v(y)|u(y)|^{2^*_\mu - 1}}{|x- y|^{n - 2}}\; \d y + u_{R/2}(y)\\
		& \leq \|v\|_{L^\infty(B_{R/2})}\int_{B_{R/2}}\frac{|u(y)|^{2^*_\mu - 1}}{|x- y|^{n- 2}}\; \d y + u_{R/2}(x), 
	\end{split}
	\end{equation}
	where $u_{R/2}$ is as in \eqref{eq:ur}. In particular, $u_{R/2}\in L^\infty_{\loc}(B_{R/2})$. Moreover, for any $x\in B_R$, H\"older's inequality and inequality \eqref{eq:case1_q_large} give
	\begin{equation}
	\label{eq:case1_simple_holders}
	\begin{split}
		\int_{B_{R/2}}\frac{|u(y)|^{2^*_\mu - 1}}{|x- y|^{n - 2}}\; \d y
		& \leq \|u\|_{L^q(B_R)}^{2^*_\mu - 1}\left(\int_{B_{2R}(x)}|x- y|^{-\frac{(n - 2)q}{q +1 - 2^*_\mu}}\; \d y\right)^{1 - \frac{2^*_\mu - 1}q}\\
		&\leq C(n, \mu, q)R^{n(\frac 2n - \frac{2^*_\mu - 1}q)}\|u\|_{L^q(B_R)}^{2^*_\mu - 1}. 
	\end{split}
	\end{equation}
	Bringing this estimate back to \eqref{eq:case1_absu_pointwise} and in view of the containment $u_{R/2}\in L^\infty_{\loc}(B_{R/2})$ we conclude that $u\in L^\infty_{\loc}(B_{R/2})$, thereby completing the proof of the lemma in \ref{item:origin_bound_case1}. 
	\item \label{item:origin_bound_case2} Assume $\mu\in (0, n - 2)$. In this case we can still choose $q$ large enough so that \eqref{eq:case1_q} is satisfied, so arguing as in the beginning of \ref{item:origin_bound_case1} we will obtain the boundedness of $v$ near the origin. However, we cannot necessarily choose $q$ to satisfy \eqref{eq:case1_q_large}, so the boundedness of $u$ near the origin does not follow from a simple application of H\"older's inequality as it did in \ref{item:origin_bound_case1}. Instead, we apply the HLS inequality in an iterative manner and improve the integrability of $u$ at each iteration. After finitely many iterations, we will have established sufficient integrability of $u$ so that the H\"older inequality argument of \ref{item:origin_bound_case1} (as in estimates \eqref{eq:case1_absu_pointwise} and \eqref{eq:case1_simple_holders}) can be applied to yield the boundedness of $u$ near the origin. The details follow. If $\mu\geq 4$ then $\sum_{k = 1}^\infty(2^*_\mu - 1)^{-k} = +\infty$ and if $\mu \in (0, 4)$ then 
	\begin{equation*}
		\frac 2n\sum_{k = 1}^\infty(2^*_\mu - 1)^{-k}
		= \frac{1}{2^*} - \frac\mu{2n} + \frac{\mu(n + 2 - \mu)}{2n(4 - \mu)}
		> \frac{1}{2^*} - \frac\mu{2n}, 
	\end{equation*}
	so
	\begin{equation*}
		N= \min\left\{m\in \bb N: \frac 2n\sum_{k = 1}^m(2^*_\mu - 1)^{-k}> \frac 1{2^*} - \frac\mu{2n}\right\}
	\end{equation*}
	is well-defined. If $N= 1$ then then we may choose $q$ for which both of estimates \eqref{eq:case1_q} and \eqref{eq:case1_q_large} hold and thus, the argument of \ref{item:origin_bound_case1} guarantees the existence of some open neighborhood of the origin on which both $u$ and $v$ are bounded. Assume that $N\geq 2$ and, for $\mu\in (0, n - 2)$ define
	\begin{equation*}
		H_\ell(\mu)
		= \begin{cases}
		(2^*_\mu - 1)^{-1} & \text{ if }\ell = 1\\
		(2^*_\mu - 1)^{-\ell} + \frac2n\sum_{k = 1}^{\ell - 1}(2^*_\mu - 1)^{-k} & \text{ if }\ell \in \{2, 3, \ldots\}. 
		\end{cases}
	\end{equation*}
	Elementary computations yield the estimate
	\begin{equation}
	\label{eq:necessary_Hl_inequality}
		\min\{H_\ell(\mu)\}_{\ell = 1}^{N - 1}
		> \frac 1{2^*} - \frac\mu {2n}
		\qquad \text{ for all }\mu\in (0, n - 2). 
	\end{equation}
	For the reader's convenience, the details of these computations are provided here.   A direct computation shows that for all $\mu\in (0, n - 2)$ and all $\ell = 1, 2, \ldots$, 
	\begin{equation*}
		H_{\ell + 1}(\mu) - H_\ell (\mu)
		= \frac{\mu - 2}{n(2^*_\mu - 1)^{\ell +1}}. 
	\end{equation*}
	Therefore, $H_{\ell + 1}(\mu)< H_\ell(\mu)$ whenever $\mu\in (0, \min\{2, n - 2\})$ and, if $n > 4$, then $H_{\ell + 1}(\mu) \geq H_\ell(\mu)$ for all $\mu\in [2, n- 2)$. If $\mu\in (0, \min\{2, n - 2\})$ then $2^*_\mu > 2$ so $H_1(\mu)> H_2(\mu)> \ldots > H_\infty(\mu)$ where 
	\begin{equation*}
		H_\infty(\mu)
		= \lim_{\ell\to\infty}H_\ell(\mu)
		= \frac{2(n - 2)}{n(4 - \mu)}. 
	\end{equation*}
	Moreover, still for $\mu\in (0, \min\{2, n- 2\})$, 
	\begin{equation*}
		H_\infty(\mu) - \left(\frac 1{2^*} - \frac{\mu}{2n}\right)
		= \frac{\mu(n + 2 - \mu)}{2n(4 -\mu)}
		> 0
	\end{equation*}
	and thus $H_\ell(\mu)> \frac{1}{2^*}- \frac{\mu}{2n}$ for all $\ell$ and all $\mu\in (0, \min\{2, n - 2\})$. If $n> 4$ and if $\mu\in [2, n- 2)$ then for any $\ell \geq 1$, 
	\begin{equation*}
	\begin{split}
		H_\ell(\mu) - \left(\frac 1{2^*} - \frac\mu{2n}\right)
		& \geq H_1(\mu) - \left(\frac 1{2^*} - \frac\mu{2n}\right)\\
		& = \frac{1}{2n(n + 2 - \mu)}\left(n^2 + 2(\mu - 2)n - (\mu - 2)(\mu + 2)\right)\\
		& \geq \frac{1}{2n(n + 2 - \mu)}\left(n^2 + 2(\mu - 2)(\mu + 2) - (\mu - 2)(\mu + 2)\right)\\
		& = \frac{n^2 + (\mu - 2)(\mu + 2)}{2n(n + 2 - \mu)}\\
		& \geq\frac{n}{2(n + 2 - \mu)}\\
		&> 0. 
	\end{split}
	\end{equation*}
	This establishes inequality \eqref{eq:necessary_Hl_inequality}. Now choose $q_1$ such that
	\begin{equation*}
		\frac 1{2^*} - \frac\mu{2n}
		< \frac 1{q_1}
		< \min\left\{\frac 2n\sum_{k = 1}^N(2^*_\mu - 1)^{-k}, \min\{H_\ell(\mu)\}_{\ell = 1}^{N - 1}, \frac{n -\mu}n\cdot\frac{1}{2^*_\mu}\right\}
	\end{equation*}	
	and define $q_2, \ldots, q_N$ by 
	\begin{equation}
	\label{eq:the_qks}
		\frac{1}{q_{k+ 1}} = \frac{2^*_\mu - 1}{q_k} - \frac 2n
	\end{equation}
	so that 
	\begin{equation}
	\label{eq:qN_large}
		\frac 1{q_N}< \frac 2n\cdot\frac{1}{2^*_\mu - 1}
	\end{equation}
	and
	\begin{equation}
	\label{eq:qk_small}
		\frac 1{q_k}> \frac 2n\cdot\frac{1}{2^*_\mu - 1}
		\qquad \text{ for }k = 1, \ldots, N - 1. 
	\end{equation}
	Observe that
	\begin{equation}
	\label{eq:implies_admissible_exponent}
		\frac 2n< \frac{2^*_\mu - 1}{q_k} < 1
		\qquad \text{ for }k = 1, \ldots, N -1
	\end{equation}
	and thus, for any such $k$, 
	\begin{equation}
	\label{eq:admissible_exponent}
		0< \frac{1}{q_{k + 1}}< 1
		\qquad \text{ and }\qquad 
		0< \frac{2^*_\mu - 1}{q_k}< 1. 
	\end{equation}
	Specifically, the first inequality in \eqref{eq:implies_admissible_exponent} follows from \eqref{eq:qk_small} and the second inequality in \eqref{eq:implies_admissible_exponent} follows from the inequality $\frac1{q_1}< \min\{H_\ell(\mu)\}_{\ell = 1}^{N - 1}$. For any positive numbers $0< r< R< \frac{3r}2$, Lemma \ref{lemma:h_integrability} guarantees that the function $h_{r, R}$ defined in \eqref{eq:h_rR} satisfies $h_{r, R}\in L^{\frac{2n}{n - 2-\mu}}(B_{R/2})\subset L^{q_1}(B_{R/2})$, so Lemma \ref{lemma:dual_integrability_boosting} guarantees the existence of $\epsilon = \epsilon(n, \mu, q_1)>0$ independent of $R$ for which $u\in L^{q_1}(B_{R/4})$ whenever $\|u\|_{L^{2^*}(B_R)}^{2^*_\mu - 1}< \epsilon$. Fix any such $\epsilon$ and choose $R = R(\epsilon)$ sufficiently small so that $\|u\|_{L^{2^*}(B_{4R})}^{2^*_\mu - 1}< \epsilon$ and hence $u\in L^{q_1}(B_R)$. For this $R$, decompose $v(x)$ as in \eqref{eq:v_Br4_decomposition}. Similarly to \ref{item:origin_bound_case1} we have $v_R\in L^\infty_{\loc}(B_R)$. Moreover, for any $x\in B_R$, similarly to \eqref{eq:case1_v-vR_pointwise}, H\"older's inequality and the inequality $\frac1{q_1}< \frac{n - \mu}n\cdot\frac{1}{2^*_\mu}$ give
	\begin{equation*}
	\begin{split}
		|v(x) - v_R(x)|
		& 
		\leq
		C(n, \mu, q_1)R^{n(\frac{n - \mu}n - \frac{2^*_\mu}{q_1})}\|u\|_{L^{q_1}(B_R)}^{2^*_\mu}, 
	\end{split}
	\end{equation*}
	so we deduce that $v\in L^\infty_{\loc}(B_R)$. Similarly to \eqref{eq:case1_absu_pointwise} we have the pointwise estimate
	\begin{equation}
	\label{eq:case2_u_pointwise_bound1}
		|u(x)|
		\leq \|v\|_{L^\infty(B_{R/2})}I_{n - 2}[|u|^{2^*_\mu - 1}\chi_{B_{R/2}}](x) + u_{R/2}(x)
	\end{equation}
	with $u_{R/2}\in L^\infty_{\loc}(B_{R/2})$. Moreover, in view of \eqref{eq:the_qks} and \eqref{eq:admissible_exponent}, the HLS inequality gives
	\begin{equation*}
		\|I_{n - 2}[|u|^{2^*_\mu - 1}\chi_{B_{R/2}}]\|_{L^{q_2}(B_{R/2})}
		\leq C\|u\|_{L^{q_1}(B_{R/2})}^{2^*_\mu - 1}
	\end{equation*}
	for some constant $C = C(n, \mu, q_1)>0$. Bringing this estimate back to \eqref{eq:case2_u_pointwise_bound1} we find that 
	\begin{equation}
	\label{eq:Lq2}
		u\in L^{q_2}_{\loc}(B_{R/2})\subset L^{q_2}(B_{R/4}).
	\end{equation}
	For $k = 2, \ldots, N$ let us define $R_k = 2^{-k}R$. Thus, in this notation, \eqref{eq:Lq2} guarantees that $u\in L^{q_2}(B_{R_2})$. For $k = 2, \ldots, N$ we have the pointwise estimate
	\begin{equation}
	\label{eq:u_kth_pointwise_estimate}
		|u(x)|
		\leq \|v\|_{L^\infty(B_{R_k})}I_{n -2}[|u|^{2^*_\mu - 1}\chi_{B_{R_k}}](x)+ u_{R_k}(x)
	\end{equation} 
	with $u_{R_k}\in L^\infty_{\loc}(B_{R_k})$. Moreover, in view of \eqref{eq:the_qks} and \eqref{eq:admissible_exponent}, for $k = 2, \ldots, N - 1$, having shown (inductively) that $u\in L^{q_k}(B_{R_k})$, the HLS inequality gives
	\begin{equation*}
		\|I_{n - 2}[|u|^{2^*_\mu - 1}\chi_{B_{R_k}}]\|_{L^{q_{k + 1}}(B_{R_{k + 1}})}
		\leq C\|u\|_{L^{q_k}(B_{R_k})}^{2^*_\mu- 1}
	\end{equation*}
	for some constant $C = C(n, \mu,q_1)>0$. Bringing this back to \eqref{eq:u_kth_pointwise_estimate} we find that $u\in L^{q_{k + 1}}(B_{R_{k + 1}})$. In particular, $u\in L^{q_N}(B_{R_N})$. Finally, for any $x\in B_{R_N}$ consider the pointwise estimate \eqref{eq:u_kth_pointwise_estimate} with $k =N$. From H\"older's inequality and \eqref{eq:qN_large}, for any $x\in B_R$ we have
	\begin{equation*}
	\begin{split}
		0
		& \leq I_{n - 2}[|u|^{2^*_\mu - 1}\chi_{B_{R_N}}](x)\\
		& \leq \|u\|_{L^{q_N}(B_{R_N})}^{2^*_\mu - 1}\left(\int_{B_{2R}(x)}|x- y|^{-\frac{q_N(n - 2)}{q_N +1 - 2^*_\mu}}\; \d y\right)^{1 - \frac{2^*_\mu - 1}{q_N}}\\
		& \leq C(n, \mu, q_1)R^{n(\frac 2n- \frac{2^*_\mu - 1}{q_N})}\|u\|_{L^{q_N}(B_{R_N})}^{2^*_\mu - 1}. 
	\end{split}
	\end{equation*}
	Bringing this estimate back to \eqref{eq:u_kth_pointwise_estimate} (still with $k = N$) and in view of the containment $u_{R_N}\in L^\infty_{\loc}(B_{R_N})$ we find that $u\in L^\infty(B_{R_{N+ 1}})$. 
\end{enumerate}
\end{proof}

\begin{proof}[Proof of Lemma \ref{lemma:local_boundedness}]
The proof follows routinely from Lemma \ref{lemma:bounded_near_origin} and the translation invariance of \eqref{eq:the_system}. 
\end{proof}
Having established the local boundedness of solutions to \eqref{eq:the_system}, the proof of Proposition \ref{prop:system_regularity} can be established by standard methods. For the reader's convenience, an outline of the proof is sketched here. 
\begin{proof}[Outline of the Proof of Proposition \ref{prop:system_regularity}]
Fix any solution $(u,v)\in L^{2^*}(\bb R^n)\times L^{\frac{2n}\mu}(\bb R^n)$ to \eqref{eq:the_system}. Lemma \ref{lemma:local_boundedness} guarantees that both of $u$ and $v$ are locally bounded on $\bb R^n$. For $R>0$ decompose $u(x)$ as
\begin{equation}
\label{eq:u_true_decomposition}
	u(x) 
	= \int_{B_R}\frac{v(y)|u(y)|^{2^*_\mu - 2}u(y)}{|x- y|^{n - 2}}\; \d y
	+ \int_{\bb R^n\setminus B_R}\frac{v(y)|u(y)|^{2^*_\mu - 2}u(y)}{|x- y|^{n - 2}}\; \d y
\end{equation}
and decompose $v(x)$ as in \eqref{eq:v_Br4_decomposition}. Routine arguments show that both the second summand on the right-hand side of \eqref{eq:u_true_decomposition} and $v_R$ are locally uniformly Lipschitz continuous on $B_R$. Moreover, by performing computations similar to those carried out in the proof of Corollary 1 of \cite{JinLiXiong2017}, one can show both that the first summand on the right-hand side of \eqref{eq:u_true_decomposition} is in $C^{0, \alpha}(\overline B_{R/4})$ whenever $\alpha\in (0, 1)$ and that $v - v_R\in C^{0, \alpha}(\overline B_{R/4})$ whenever $\alpha\in (0, \min \{1, n - \mu\})$. Therefore, $u\in C^{0, \alpha}(\overline B_{R/4})$ whenever $\alpha\in (0, 1)$ and $v\in C^{0, \alpha}(\overline B_{R/4})$ whenever $\alpha\in (0, \min\{1, n - \mu\})$. Using these containments and using the fact that $R>0$ is arbitrary we find that for any H\"older exponent $\alpha\in (0, \min\{1, n - \mu\})$ the function 
\begin{equation}
\label{eq:integrand_f}
	f = v|u|^{2^*- 2}u
\end{equation}
is locally uniformly H\"older continuous of exponent $\alpha$ on $\bb R^n$. Furthermore, the integrability assumptions on $u$ and $v$ guarantee that $f\in L^{\frac{2n}{n + 2}}(\bb R^n)$. Now let $\eta\in C^\infty([0, \infty);[0, 1])$ satisfy both $\eta\equiv 0$ on $[0, 1]$ and $\eta\equiv 1$ on $[2, \infty)$ and for $R>0$ define $\eta_R\in C^\infty(\bb R^n; [0, 1])$ by $\eta_R(\xi) = \eta(R^{-1}|\xi|)$. For $f$ as in \eqref{eq:integrand_f} we have $u(x) = J_1(x) + J_2(x)$, where
\begin{equation*}
\begin{split}
	J_1(x) 
	& = \int_{\bb R^n}\frac{\eta_R(x- y)}{|x- y|^{n - 2}}f(y)\; \d y\\
	J_2(x) 
	& = \int_{\bb R^n}\frac{1 - \eta_R(x- y)}{|x- y|^{n - 2}}f(y)\; \d y. 
\end{split}
\end{equation*}
For any $\alpha\in (0, \min\{1, n - \mu\})$, using the fact that $f$ is uniformly H\"older continuous of exponent $\alpha$ on $\overline B_{2R + 2}$, the existence and uniform H\"older continuity of $D^2J_2$ with exponent $\alpha$ on $\overline B_{R/2}$ can be established using the techniques of the proofs of Lemmata 4.1-4.2 of \cite{GilbargTrudinger1977}. Using the fact that the integral kernel $K(\xi) = \eta_R(\xi)|\xi|^{2 - n}$ in the definition of $J_1$ is smooth and vanishes identically on $B_R$ together with the containment $f\in L^{\frac{2n}{n + 2}}(\bb R^n)$ one can show that $D^2J_1$ is uniformly Lipschitz continuous on $\overline B_{R/2}$. 
\end{proof}
\bibliographystyle{alpha}


\end{document}

%% file: structure.tex
\usepackage[T1]{fontenc} 
\usepackage[utf8]{inputenc} 
\usepackage{amssymb}
\usepackage{graphicx} 
\graphicspath{{Figures/}} 
\usepackage{enumitem} 
\usepackage{subfig} 
\usepackage{varioref} 
\usepackage[authoryear,sort&compress]{natbib}
\usepackage{tikz}
\usepackage{tikz-cd}
\usetikzlibrary{arrows,arrows.meta}
\usetikzlibrary{calc}
\usepackage{esint} 
\usepackage{amsaddr} 
\theoremstyle{plain} 
\newtheorem{theorem}{Theorem}
\newtheorem{lemma}[theorem]{Lemma}
\newtheorem{prop}[theorem]{Proposition}
\newtheorem{coro}[theorem]{Corollary}

\newtheorem{oldtheorem}{Theorem}

\theoremstyle{definition} 
\newtheorem{defn}[theorem]{Definition}
\newtheorem{example}[theorem]{Example}

\theoremstyle{remark} 
\newtheorem{remark}[theorem]{Remark}

\numberwithin{equation}{section}
\numberwithin{theorem}{section}
\usepackage[bookmarks, colorlinks=true, pdfstartview=FitV, linkcolor=blue, citecolor=blue, urlcolor=blue]{hyperref}

\renewcommand{\d}{\mathrm d}
\renewcommand{\div}{\mathrm{div}}

\newcommand{\lap}{\Delta}
\newcommand{\Grad}{\nabla}
\newcommand{\weakconv}{\rightharpoonup}

\newcommand{\mc}{\mathcal}
\newcommand{\bb}{\mathbb}
\newcommand{\abs}[1]{\left|#1\right|}
\newcommand{\norm}[1]{\left\|#1\right\|}
\newcommand{\lb}{\langle}
\newcommand{\rb}{\rangle}

\DeclareMathOperator{\card}{card}

\DeclareMathOperator{\Id}{Id}

\DeclareMathOperator{\loc}{loc}
